\documentclass[seceqn,secthm]{elsart}
\usepackage{amssymb}
\usepackage{amsfonts}
\usepackage{amsmath}
\usepackage{graphicx}%
\numberwithin{equation}{section}
\begin{document}

\begin{frontmatter}
\title{Approximation of Quantum Tori by Finite Quantum Tori for the quantum
Gromov-Hausdorff distance}
\author{Fr\'{e}d\'{e}ric Latr\'{e}moli\`{e}re}
\address{Department of Mathematics, University of Toronto}
\ead{frederic@math.toronto.edu}

\thanks{The author's research was supported in part by NSF Grant DMS-0200591.}

\date{October, 15$^{\text{th}}$ 2004.}

\begin{abstract}
We establish that, given a compact Abelian group $G$ endowed with a continuous
length function $l$ and a sequence $(H_{n})_{n\in\mathbb{N}}$ of closed
subgroups of $G$ converging to $G$ for the Hausdorff distance induced by $l$,
then $C^{\ast}\left(  \widehat{G},\sigma\right)  $ is the quantum
Gromov-Hausdorff limit of any sequence $C^{\ast}\left(  \widehat{H_{n}}%
,\sigma_{n}\right)  _{n\in\mathbb{N}}$ for the natural quantum metric
structures and when the lifts of $\sigma_{n}$ to $\widehat{G}$ converge
pointwise to $\sigma$. This allows us in particular to approximate the quantum
tori by finite dimensional C*-algebras for the quantum Gromov-Hausdorff
distance. Moreover, We also establish that if the length function $l$ is
allowed to vary, we can collapse quantum metric spaces to various quotient
quantum metric spaces.

\end{abstract}


\end{frontmatter}

\section{Introduction}

The quantum Gromov-Hausdorff distance defines a framework to investigate
claims found in the physics literature that classical or quantum spaces can be
approximated by Matrix algebras, as in \cite{Connes97}, \cite{Madore} for
instance. These approximations are related to M-theory \cite{Connes97},
\cite{Zumino98}, \cite{Szabo01} and may serve as tools to construct quantum
field theory over quantum spaces, generalizing methods found in \cite{Jaffe}.
In this paper, we specifically investigate the matter of approximating the
quantum tori by matrix algebras generated by unitaries in finite dimension,
which is suggested for instance in \cite[sec. 3.3]{Connes97}.

We start by recalling the foundations of the quantum Gromov-Haus\-dorff
distance theory. Rieffel introduced in \cite{Rieffel00} a notion of
convergence for compact quantum metric spaces \cite{Rieffel99}, which
generalizes to noncommutative geometry the Gromov-Hausdorff distance
\cite[Ch.3]{Gromov} between isometry classes of compact metric spaces. A
compact quantum metric space $(A,L)$, as defined in \cite{Rieffel99} and
\cite{Rieffel00}, is an order-unit space $A$ endowed with a seminorm $L\;$such
that the dual distance $d_{L}$, defined for any two states $\mu,\nu$ of $A$
by
\begin{equation}
d_{L}(\mu,\nu)=\sup\{\left\vert \mu(a)-\nu(a)\right\vert :a\in A\text{,
}L(a)\leq1\}\text{,}\label{dual_metric}%
\end{equation}
induces the weak* topology on the state space $S(A)$ of $A$. Such a seminorm
is called a Lip-norm. The classical examples of Lip-norms are the Lipschitz
seminorms on $C(X)\ $for any compact metric space $(X,d)$ (see \cite[sec.
11.3, Theorem 11.3.3]{Dudley}). Rieffel's definition refines an idea in
\cite{Connes89}\cite[Chapter VI]{Connes}.

Now, let $(A,L)$ and $(B,L^{\prime})$ be two compact quantum metric spaces
with respective state spaces $S(A)\ $and $S(B)$. Let $\mathcal{C}(L,L^{\prime
})\ $be the set of Lip-norms on the order-unit space $A\oplus B$ whose
quotients on $A$ (resp. $B$) are the Lip-norm $L$ (resp. $L^{\prime}$).
Rieffel proves in \cite{Rieffel00} that $\mathcal{C}(L,L^{\prime})$ is
nonempty, and defines the quantum Gromov-Hausdorff distance by setting:%
\[
\operatorname{dist}_{q}((A,L),(B,L^{\prime}))=\inf\left\{  \mathfrak{H}%
[d_{L^{\prime\prime}}]\ (S(A),S(B)):L^{\prime\prime}\in\mathcal{C}%
(L,L^{\prime})\right\}  \text{,}%
\]
where $\mathfrak{H}[d_{L^{\prime\prime}}]$ is the Hausdorff distance induced
by the distance $d_{L^{\prime\prime}}$ on the weak*-closed subsets of the
state space $S(A\oplus B)$ of $A\oplus B$, and where $d_{L^{\prime\prime}}$ is
defined by (\ref{dual_metric}).

We now turn to the problem resolved in this paper. Let $G$ be a compact
Abelian group endowed with a continuous length function $l$, and denote by $e$
the unit of $G$. The length function $l$ induced a metric $d_{l}$ on the group
$G$ by setting $d_{l}(g^{\prime},g)=l(g^{-1}g^{\prime})$ for all $g,g^{\prime
}\in G$, and $d_{l}$ in turn generates a topology $\tau_{2}$ on $G$. We denote
by $\tau_{1}$ the original (compact) topology of $G$, and we denote by
$\mathfrak{H}(l)$ the Hausdorff distance defined by $d_{l}$ on the closed
subsets of $G$. We make the following simple observation:

\begin{lem}
\label{dense}We have: $\tau_{1}=\tau_{2}$. The group $G$ is therefore
separable. In particular, if $(X_{n})_{n\in\mathbb{N}}$ is a sequence of
closed subsets of $G$ converging to $G$ for $\mathfrak{H}(l)$, then $%
{\displaystyle\bigcup\limits_{n=0}^{\infty}}
X_{n}$ is dense in $G$.
\end{lem}

\begin{pf}
Since $l$ is continuous, it is immediate that $\tau_{2}\subseteq\tau_{1}$. The
other inclusion is from a standard argument \cite{Rudin91} which we include
here for convenience. Let $F$ be a closed set for $\tau_{1}$. Since
$(G,\tau_{1})$ is compact, $F$ is compact as well. Since $\tau_{2}%
\subseteq\tau_{1}$, the set $F$ is compact for $\tau_{2}$ as well. Since
$\tau_{2}$ is Hausdorff, the set $F$ is closed in $(G,\tau_{2})$. Hence
$\tau_{1}=\tau_{2}$. The topology $\tau_{1}$ is compact and metrizable under
our assumptions, hence $G$ is separable. Now, since $(X_{n})_{n\in\mathbb{N}}$
converges to $G$ for $\mathfrak{H}(l)$, we have by definition that for any
$\varepsilon>0$ and $g\in G$ there exists $N\in\mathbb{N}$ such that for all
$n\geq N$ the $\tau_{2}$-open ball $\{x\in G:l(x^{-1}g)<\varepsilon\}$
intersects $X_{n}$. Hence all open balls in $\tau_{2}$ intersect $%
{\displaystyle\bigcup\limits_{n=0}^{\infty}}
X_{n}$, which is therefore dense in $\tau_{2}$, hence in $\tau_{1}$. \qed \end{pf}

\begin{rem}
Lemma (\ref{dense}) shows that we do not lose any generality by working with a
group $G$ endowed with a length function $l$ which generates a compact
topology on $G$, or with a separable compact group $G$, which can always be
endowed with a continuous length function. Varying the length function on $G$
is the subject of Theorem (\ref{Collapse}).
\end{rem}

The Pontryagin topological dual group of $G$ is denoted by $\widehat{G}.$ Of
course, $\widehat{G}$ is discrete. We denote the group of skew bicharacters of
$\widehat{G}$ by $A(\widehat{G})$, and endow it with the topology of pointwise
convergence over $\widehat{G}$. The group $A(\widehat{G})$ is compact by
\cite[Prop. 2.1 and Theorem 2.1]{Kleppner65}. Let $\sigma_{\infty}$ be a
skew-bicharacter of $\widehat{G}$.

Let now $(H_{n})_{n\in\mathbb{N}}$ be a sequence of closed subgroups of $G$,
converging to $G$ for $\mathfrak{H}(l)$. For each $n\in\mathbb{N}$, we denote
the dual group of $H_{n}$ by $\widehat{H_{n}}$, which is a quotient group of
$\widehat{G}$ and is discrete (again as the dual of a compact group). We
denote by $J_{n}$ the closed subgroup of $\widehat{G}$ such that
$\widehat{H_{n}}=\widehat{G}/J_{n}$. Equivalently, $J_{n}$ is the group of
characters of $G$ which are identically 1 on $H_{n}$ (the annihilator of
$H_{n}$).

To ease notations, we will denote by $H_{\infty}=G,$ so that $\widehat
{H_{\infty}}=\widehat{G}$, and accordingly we set $J_{\infty}=\{1\}$, where
$1$ is the trivial character of $G$. We will denote by $\overline{\mathbb{N}%
}=\mathbb{N}\cup\{\infty\}$ the one-point compactification of $\mathbb{N}.$

For each $n\in\overline{\mathbb{N}}$, we denote by $\sigma_{n}$ a
skew-bicharacter on $\widehat{H_{n}}$. This bicharacter naturally lifts to
$\widehat{G}$ to a skew bicharacter $\sigma_{n}^{\prime}$ on $\widehat{G}$. We assume
that $(\sigma_{n}^{\prime})_{n\in\mathbb{N}}$ converges to $\sigma_{\infty}$
in $A(\widehat{G})$.

Let $n\in\overline{\mathbb{N}}$. The canonical dual action $\alpha_{H_{n}}$ of
$H_{n}$ on $C^{\ast}\left(  \widehat{H_{n}},\sigma_{n}\right)  $ is an ergodic
strongly continuous action obtained by extending by continuity the map:%
\[
\alpha_{H_{n}}^{g}:f\in C_{c}\left(  \widehat{H_{n}}\right)  \longmapsto
\left\langle g,.\right\rangle f\in C_{c}\left(  \widehat{H_{n}}\right)
\]
for all $g\in H_{n}$, where $C_{c}(X)$ is the space of compactly supported
continuous functions on any locally compact space $X$ and $\left\langle
.,.\right\rangle $ is the dual pairing between $H_{n}$ and $\widehat{H_{n}}$.
Observe that on $C_{c}(\widehat{H_{n}})$, the restriction of $\alpha_{G}$ to
$H_{n}$ is $\alpha_{H_{n}}$. We denote by $\left\Vert .\right\Vert _{n}$ the
C*-norm of $C^{\ast}\left(  \widehat{H_{n}},\sigma_{n}\right)  $ (see
\cite{Zeller-Meier68}). By \cite{Rieffel98a}, the seminorm defined for all
$a\in C^{\ast}\left(  \widehat{H_{n}},\sigma_{\infty}\right)  $ by%
\[
L_{n}(a)=\sup\left\{  l(g)^{-1}\left\Vert a-\alpha_{H_{n}}^{g}(a)\right\Vert
_{n}:g\in H_{n}\backslash\{e\}\right\}
\]
is a Lip-norm on $C^{\ast}\left(  \widehat{H_{n}},\sigma_{n}\right)
^{\operatorname*{sa}}$, where for any subset $X$ of $C^{\ast}\left(
\widehat{H_{n}},\sigma_{n}\right)  $ we denote by $X^{\operatorname*{sa}}$ the
set of selfadjoint elements in $X$.

We shall prove the following theorem under these assumptions:

\begin{thm}The sequence $\left( C^{\ast }\left(
\widehat{H_{n}},\sigma_{n}\right)  ^{\operatorname*{sa}} ,L_{n}\right)  
_{n\in\mathbb{N}}$ of compact quantum metric spaces converges to
$\left( C^{\ast}\left( \widehat{G},\sigma_{\infty}\right)  
^{\operatorname*{sa}},L_{\infty}\right)  $ for the quantum 
Gromov-Hausdorff distance. \end{thm}

Our result is parallel to the approximation by full matrix algebras of $C(O)$
for any coadjoint orbit $O$ of any compact semi-simple Lie group $G$ developed
by Rieffel in \cite{Rieffel01}. Note also that the problem solved in this
paper was raised in \cite{Rieffel01}.

We recall from \cite{Zeller-Meier68} that if $\theta$ and $\theta^{\prime}$
are cohomologous 2-cocycles of our discrete group $\widehat{G}$, then the
C*--algebras $C^{\ast}(\widehat{G},\theta)$ and $C^{\ast}(\widehat{G}%
,\theta^{\prime})$ are *-isomorphic. Moreover, as $\widehat{G}$ is Abelian,
discrete and separable, given any 2-cocycle $\theta$ of $\widehat{G}$ we can
find a skew bicharacter $\sigma$ of $\widehat{G}$ such that $\sigma$ is
cohomologous to $\theta$ by \cite[Theorem 7.1]{Kleppner65}. Up to
*-isomorphism, all the C*--algebras $C^{\ast}(\widehat{G},\theta)$, for
$\theta$ a 2-cocycle of $\widehat{G}$, are thus of the form $C^{\ast}%
(\widehat{G},\sigma)$, where $\sigma$ is a skew bicharacter of $\widehat{G}$.
It is easier to work with convergence of bicharacters than the less intuitive
notion of convergence of cocycles, which justifies our particular choice for
representing 2-cocycles.

The fundamental example we have in mind for this situation is the following:

\begin{exmp}
[Quantum Torus]\label{fundamentalexample}Let $G=\mathbb{T}^{d}$ be the
$d$-dimensional torus for some $d\in\{1,2,\ldots\}$. We suppose given a
sequence $(k_{n})_{n\in\mathbb{N}}$ of elements in $\mathbb{N}^{d}$ such that
$\lim_{n\rightarrow\infty}k_{n}=(\infty,\ldots,\infty)$. For any
$k\in\mathbb{N}^{d}$, written as $k=(k(1),\ldots,k(d))$, we define
$k\mathbb{Z}^{d}=%
{\displaystyle\prod\limits_{j=1}^{d}}
k(j)\mathbb{Z}$ and $\mathbb{Z}_{k}^{d}=\mathbb{Z}^{d}/k\mathbb{Z}^{d}$. Note
that up to isomorphism, all quotient groups of $\mathbb{Z}^{d}$ are obtained
this way. Now, we set $\widehat{H_{n}}=\mathbb{Z}_{k_{n}}^{d}$ for all
$n\in\mathbb{N}$, and we denote $H_{n}$ by $\mathbb{U}_{k_{n}}^{d}$. Let
$S_{\infty}$ be $d\times d$ antisymmetric matrix and for all $n\in\mathbb{N}$
let $S_{n}=((a_{i,j,n}))_{1\leq i,j\leq d}$ be an antisymmetric matrix such
that $a_{i,j,n}\gcd(k(i),k(j))\in\mathbb{Z}$, with the convention that
$\gcd(0,m)=\gcd(m,0)=m$ for any $m\in\mathbb{N}$. Let $\cdot$ be the canonical
dot product on $\mathbb{R}^{d}$. It is easy to check that $\left(  \chi
,\chi^{\prime}\right)  \in\widehat{H_{n}}\times\widehat{H_{n}}\mapsto
\exp(2i\pi S_{n}\chi\cdot\chi^{\prime})$ is a skew bicharacter $\sigma\lbrack
S_{n}]$ of $\widehat{H_{n}}$ for all $n\in\overline{\mathbb{N}}$ when we embed
$\widehat{H_{n}}$ naturally in $\mathbb{R}^{d}$. Assume at last that
$(S_{n})_{n\in\mathbb{N}}$ converges in operator norm to $S_{\infty}$. Our
main result in this paper implies that, for any continuous length function on
$\mathbb{T}^{d}$, the "fuzzy"\ tori $C^{\ast}\left(  \mathbb{Z}_{k_{n}}%
^{d},\sigma\lbrack S_{n}]\right)  _{n\in\mathbb{N}}$ converge to the quantum
torus $C^{\ast}\left(  \mathbb{Z}^{d},\sigma\lbrack S_{\infty}]\right)  $ when
$n\rightarrow\infty$.
\end{exmp}

To prove Theorem (\ref{MAIN}), we will first prove that the family of C*-algebras
$C^{\ast}(\widehat{H_{n}},\sigma_{n})_{n\in\overline{\mathbb{N}}}$ admit a
continuous field structure using groupoid techniques. We then prove our
theorem with some harmonic analysis. In the last section of this paper, we
address the natural question of convergence when the length function $l$ is
allowed to vary. In particular, we show that it is possible to prove that one
can collapse a quantum torus to a lower dimensional quantum torus. The results
of our last section hold for the more general setting of a compact group $G$
acting ergodically and strongly continuously on any unital C*-algebra.

\section{Continuous fields of twisted groupoids C*-algebras}

For each $n\in\mathbb{N}$, The group $A_{n}$ of skew-bicharacters of
$\widehat{H_{n}}$ is identified with the closed subgroup of $A(\widehat{G})$
of the lifts to $\widehat{G}$ of the elements in $A_{n}$. By \cite[Proposition
2.2]{Kleppner65}, the compact group $A(\widehat{G})$ is a metrizable space.

\subsection{Construction of a groupoid}

We will use the notations and terminology of \cite{Renault80}. Let us start by
defining a trivial group bundle $\Omega=\Omega^{(0)}\times\widehat{G}$ with
its natural groupoid structure \cite{Renault80}, where the space of
units$\;\Omega^{(0)}$ is defined by:%
\[
\Omega^{(0)}=\left\{  (n,\sigma)\in\overline{\mathbb{N}}\mathbb{\times
}A\left(  \widehat{G}\right)  :\sigma\in A_{n}\right\}
\]
and thus the multiplication of $(n,\sigma,\chi)\in\Omega$ and $(n^{\prime
},\sigma^{\prime},\chi^{\prime})\in\Omega$ is defined if and only if
$n=n^{\prime}$ and $\sigma=\sigma^{\prime}$ and is then given by
$(n,\sigma,\chi\chi^{\prime})$.

The groupoid $\Omega$ is topologized simply as a subset of $\overline
{\mathbb{N}}\times A(\widehat{G})\times\widehat{G}$, thus it is a metrizable
locally compact space by construction. It is straightforward that the
multiplication of the groupoid $\Omega$ is continuous on the sets of
composable pairs by continuity of the multiplication in $\widehat{G}$, and
just as easily one checks that the inverse map of $\Omega$ is continuous.
Hence, $\Omega$ is a topological groupoid as defined in \cite[Definition 2.1,
p. 16]{Renault80}. In addition, if $\mu$ is the counting measure of
$\widehat{G}$, then the constant function $u\in\Omega^{(0)}\mapsto\mu$ is
obviously a (left) Haar system on $\Omega$. Yet, we will not use this Haar
system here.

Now, the groupoid $\Omega$ contains the following subgroupoid, which is itself
a group bundle over $\Omega^{(0)}$:

\begin{lem}
$\widetilde{\Omega}=\left\{  (n,\sigma,\chi)\in\Omega:\chi\in J_{n}\right\}  $
is a closed subgroupoid of $\Omega$ and $\widetilde{\Omega}^{(0)}=\Omega
^{(0)}$.
\end{lem}

\begin{pf}
It is trivial that $\widetilde{\Omega}$ is a subgroupoid of $\Omega$ with the
same unit space $\Omega^{(0)}$, so we shall only show that $\widetilde{\Omega
}$ is closed in $\Omega$. Since $\Omega$ is metrizable, it is enough to check
that $\widetilde{\Omega}$ is sequentially closed in $\Omega$. Let
$(k_{n},\sigma_{n},\chi_{n})_{n\in\mathbb{N}}$ be a sequence in $\widetilde
{\Omega}$ which converges in $\Omega$ to $(k,\sigma,\chi)$. If $k\in
\mathbb{N}$, then there exists $N\in\mathbb{N}$ such that for all $n\geq N$ we
have $k_{n}=k$, and therefore $\sigma_{n}\in A_{k}$ and $\chi_{n}\in J_{k}$.
The group $J_{k}$ is discrete and the group $A_{k}$ is closed in
$A(\widehat{G})$, so $(\sigma,\chi)\in A_{k}\times J_{k}$ and thus
$(k,\sigma,\chi)\in\widetilde{\Omega}$.

The alternative case is $k=\infty$. Since $\widehat{G}$ is discrete, there
exists $N\in\mathbb{N}$ such that for all $n\geq N$ we have $\chi_{n}=\chi$.
Therefore, by definition of $\widetilde{\Omega}$, the character $\chi$ of $G$
is 1 on $H_{k_{n}}$ for all $n\geq N$. Yet, by assumption, $(H_{k_{n}}%
)_{n\in\mathbb{N}}$ converges to $G$ for $\mathfrak{H}(l)$, since
$(H_{n})_{n\in\mathbb{N}}$ does and $\lim_{n\rightarrow\infty}k_{n}=\infty$.
Hence by Lemma (\ref{dense}), the subset $\bigcup_{n=N}^{\infty}H_{k_{n}}$ of
$G$ is dense in $G$. Since $\chi$ is continuous on $G$ and constant equal to 1
on a dense subset of $G$, it is the trivial character 1. Hence $(k,\sigma
,\chi)=(\infty,\sigma,1)\in\widetilde{\Omega}$. \qed \end{pf}

\begin{defn}
Let $\Gamma$ be the algebraic quotient space $\widetilde{\Omega}%
\backslash\Omega$. Formally, $\Gamma$ is the quotient space of the equivalence
relation $\sim$ on $\Omega$ defined for all $x,y\in\Omega$ by%
\[
x\sim y\text{ if and only if }(x,y)\in\Omega^{(2)}\text{ and }y^{-1}%
x\in\widetilde{\Omega}\text{.}%
\]

\end{defn}

\begin{lem}
\label{main_groupoid} Endowed with the quotient topology and the quotient
algebraic structure, $\Gamma$ is an amenable separable Hausdorff locally
compact groupoid such that its space of units is $\Gamma^{(0)}=\Omega^{(0)}$,
and the canonical surjection $q:\Omega\rightarrow\Gamma$ is an open epimorphism.
\end{lem}

\begin{pf}
By \cite[proposition (2.1), p. 75]{Renault80}, $\Gamma$ is a locally compact
Hausdorff space and $q$ is open. Standard algebraic manipulations show that
$q$ defines a groupoid structure on $\Gamma$, which is again a trivial bundle
of groups over $\Omega^{(0)}$. The fiber of range (and source) $(n,\sigma
)\in\Gamma^{(0)}(=\Omega^{(0)})$ is just the discrete group $\widehat{H_{n}}$,
hence the groupoid $\Gamma$ is $r$-discrete. The map $(n,\sigma)\in
\Gamma^{(0)}\mapsto\mu_{n}$, where $\mu_{n}$ is the counting measure on
$\widehat{H_{n}}$, defines a left Haar system on $\Gamma$. We fix this system
in the rest of this paper. We remark that this Haar system is not obtained as
the image of the Haar system on $\Omega$ by $q$. Last, since $G$ is separable,
so is $\widehat{G}$, so $\Omega$ is separable and so is $\Gamma$. Since
$\widehat{H_{n}}$ is Abelian for any $n\in\overline{\mathbb{N}}$, it is
amenable and thus so is $\Gamma$. \qed \end{pf}

More generally, groupoids of this type were used in \cite{Higson01}\ as
sources of counter-examples for the Baum-Connes conjecture, but this is not
the direction we will follow now.

In the case of the two dimensional quantum tori, as well as some other special
cases, it is interesting to note that a different groupoid can be constructed
and serve the same purpose of obtaining continuous fields of C*-algebras. In
fact, we will provide here a construction which is more general than what we
really need for the purpose of this paper, but is worthy of interest none the less.

We consider a sequence $(X_{n})_{n\in\mathbb{N}}$ of closed subsets of a
compact metric space $(X_{\infty},d)$, such that $(X_{n})_{n\in\mathbb{N}}$
converges to $X_{\infty}$ for the Hausdorff distance $\mathfrak{H}[d]$ defined
by $d$. For each $n\in\overline{\mathbb{N}}$, we assume given a continuous
action $\alpha_{n}$ of the group $\widehat{H_{n}}$ on the space $X_{n}$. We
assume that for any $\varepsilon>0$, any $g\in\widehat{G}$ and any $y\in
X_{\infty}$, there exists $\delta>0$ such that for all $n\in\overline
{\mathbb{N}}$, for all $x\in X_{n}$ such that $d(x,y)<\delta$ and for all
$h\in\widehat{H_{n}}$ such that $l(h^{-1}g)<\delta$, we have $d(\alpha_{n}%
^{h}(x),\alpha_{\infty}^{g}(y))\leq\varepsilon$.

Now, the set $\Theta=\left\{  (n,x,\chi):n\in\overline{\mathbb{N}}\text{,
}x\in X_{n}\text{, }\chi\in\widehat{G}\right\}  $ is a groupoid for the
multiplication defined for $(n,x,\chi)\in\Theta$ and $(n^{\prime},x^{\prime
},\chi^{\prime})\in\Theta$ only when $n=n^{\prime}$ and $x^{\prime}=\alpha
_{n}^{\chi}(x)$ by setting $(n,x,\chi)(n,\alpha_{n}^{\chi}(x),\chi^{\prime
})=(n,x,\chi\chi^{\prime})$.

In particular, the unit space $\Theta^{(0)}$ of $\Theta$ is $\left\{
(n,x,1):n\in\mathbb{N}\text{, }x\in X_{n}\right\}  $, the source map is given
by $s_{\Theta}(n,x,\chi)=(n,x,1)$ and the range map by $r_{\Theta
}(n,x,\chi)=(n,\alpha_{n}^{\chi}(x),1)$ for all $(n,x,\chi)\in\Theta$. The
inverse of $(n,x,\chi)$ is $(n,\alpha_{n}^{\chi}(x),\chi^{-1})$.

Now, $\Theta$ is topologized naturally as a subset of $\overline{\mathbb{N}%
}\times X\times\widehat{G}$. Our assumption on the coherence of the actions
$\alpha_{n}$ ensure that the multiplication in $\Theta$ and the inverse map
are continuous. Now, as before, we define two new groupoids $\widetilde
{\Theta}$ and $\widetilde{\Gamma}$ by $\widetilde{\Theta}=\left\{
(n,x,\chi)\in\Theta:\chi\in J_{n}\right\}  $ and  $\widetilde{\Gamma
}=\widetilde{\Theta}\backslash\Theta$. Again, the unit space of $\widetilde
{\Gamma}$ is $\Theta^{(0)}$. Let us define the map $p:\widetilde{\Gamma
}\rightarrow\overline{\mathbb{N}}$ simply by $p(n,x,\chi)=n$. It is immediate
that $p^{-1}(\{n\})$ is the transformation groupoid $\widehat{H_{n}}%
\times_{\alpha_{n}}X_{n}$ of the action $\alpha_{n}$ of $\widehat{H_{n}}$ on
$X_{n}$. This allows us to construct a natural left Haar system on
$\widetilde{\Gamma}$, and we will now see how this structure helps obtain a
continuous field of C*-algebras.

\begin{exmp}
[Quantum\ Tori of dimension 2]With the notations of Example (\ref{fundamentalexample}
), we set $X_{n}=\mathbb{U}_{n}^{1}$ and $X_{\infty}=\mathbb{T}$. Let
$\theta_{n}\in\mathbb{U}_{n}^{1}$ for all $n\in\mathbb{N}$ and $\theta
\in\mathbb{T}$. Then $\widehat{H_{n}}=\mathbb{Z}_{n}^{1}$ acts on $X_{n}$ by
rotation by $\theta_{n}$, and $\mathbb{Z}$ acts on $\mathbb{T}$ by rotation by
$\theta$ as well. Of course $\mathbb{Z}\times_{\theta}C(\mathbb{T})$ is
a quantum torus.
\end{exmp}

\subsection{Continuous fields of twisted groupoid C$^{\text{*}}$--algebras}

We recall the following definition from \cite[Definition 2.1.1]{Ram98}, also
in \cite[Definition 5.2]{LanRam00}:

\begin{defn}
\label{contfieldgroupoid}A \emph{continuous field of groupoids }%
$(\mathcal{G},p,T)$ is a locally compact groupoid $\mathcal{G}$ together with
a Hausdorff locally compact space $T$ and a continuous open map $p:\mathcal{G}%
\rightarrow T$ such that $p=p\circ r=p\circ s$.
\end{defn}

As a consequence of this definition, for any $t\in T$, the set $p^{-1}(\{t\})$
is a closed subgroupoid of $\mathcal{G}$, namely the reduction of
$\mathcal{G}$ to the saturated closed subset $p^{-1}(\{t\})\cap\mathcal{G}%
^{(0)}$, as defined in \cite[Definition 1.4, p.8]{Renault80}. We shall denote
this groupoid by $\mathcal{G}[t]$. It is locally compact, and it inherits by
restriction a left Haar system from any left Haar system on $\mathcal{G}$. In
our previous paragraph, we encountered two such continuous fields of
groupoids:\ $(\Gamma,r,\Omega^{(0)})$ and $(\widetilde{\Gamma},p,\overline
{\mathbb{N}})$.

Ramazan proves in \cite[Theorem 2.4.6]{Ram98}\ (see also \cite[Theorem
5.5]{LanRam00}) that $C^{\ast}(\mathcal{G})$ is then an algebra of continuous
sections for the family $(C^{\ast}(\mathcal{G}[t])_{t\in T}$. The proof of
this result essentially follows the layout proposed by \cite{Rieffel89} to
prove lower and upper semicontinuity for fields of C*--algebras. The key
result used for both those steps in the case of groupoids is the
disintegration of representations of twisted groupoid crossed products proven
by Renault in \cite[Theorem 4.1]{Renault87}.

\bigskip We will need in our work a slightly more general result:

\begin{thm}
\label{cont}Let $(\mathcal{G},p,T)$ be a continuous field of groupoid, and
assume $\mathcal{G}$ is separable. Let $\theta$ be a continuous 2-cocycle on
$\mathcal{G}$. Let $f\in C_{c}(\mathcal{G})$, and for all $t\in T$ denote by
$f^{t}$ the restriction of $f$ to the subgroupoid $\mathcal{G}[t]:=p^{-1}%
(\{t\})$ and by $\theta_{t}$ the restriction of $\theta$ to $\mathcal{G}[t]$.
For any C*-algebra $A$, we denote the norm of $A$ by $\left\Vert .\right\Vert
_{A}$.%
\begin{align*}
t  &  \longmapsto\left\Vert f^{t}\right\Vert _{C^{\ast}(\mathcal{G}%
[t],\theta_{t})}\text{ is upper semicontinuous over }T\text{,}\\
t  &  \longmapsto\left\Vert f^{t}\right\Vert _{C_{\operatorname{red}}^{\ast
}(\mathcal{G}[t],\theta_{t})}\text{ is lower semicontinuous over }T\text{.}%
\end{align*}

\end{thm}

In particular, when $\mathcal{G}$ is amenable we conclude that $t\mapsto
\left\Vert f^{t}\right\Vert _{C^{\ast}(\mathcal{G}[t],\theta_{t})}$ is
continuous over $T$, hence:

\begin{cor}
\label{cont0}$C^{\ast}(\mathcal{G},\theta)$ provides a continuous structure
for the field of C*--algebras $(C^{\ast}(\mathcal{G}[t],\theta_{t}))_{t\in T}$
when $\mathcal{G}$ is an amenable locally compact groupoid.
\end{cor}

The complete proof of these results can be found in \cite{Latremoliere04}, but
since it follows essentially well-known techniques already used in
\cite[Theorem 2.4.6]{Ram98}, we shall omit it here. The remark which allows us
to generalize \cite[Theorem 2.4.6]{Ram98} to our situation is simply that, if
$U$ is a saturated subset of $\mathcal{G}^{(0)}$, then the restricted groupoid
$\mathcal{G}_{U}$ satisfy the property that if $x\in\mathcal{G}_{U}$ and
$y\in\mathcal{G}$ then $(x,y)\in\mathcal{G}^{(2)}$ or $(y,x)\in\mathcal{G}%
^{(2)}$ if, and only if $y\in\mathcal{G}_{U}$ and either $(x,y)$ or $(y,x)$
are in $\mathcal{G}_{U}^{(2)}$. Therefore, the restriction of the continuous
2-cocycle $\theta$ to $\mathcal{G}_{U}$ is a continuous 2-cocycle of
$\mathcal{G}_{U}$. From this, the argument in \cite{Ram98} works just as well
in our setting to prove Theorem (\ref{cont}). The key observation is the
following exact sequence:

\begin{prop}
\label{exact}Let $\mathcal{G}$ be a separable locally compact groupoid. Let
$U$ be a saturated open subset of $\mathcal{G}^{(0)}$ and let $\theta$ be a
continuous 2-cocycle of $\mathcal{G}$. Let $F$ be the complement of $U$ in
$\mathcal{G}^{(0)}$, and let $\mathcal{G}_{U}$ (resp. $\mathcal{G}_{F}$) be
the reduced groupoid from $\mathcal{G}$ to $U\ $(resp. $F$). Then the
following sequence is exact:%
\[
0\longrightarrow C^{\ast}(\mathcal{G}_{U},\theta_{U})\longrightarrow C^{\ast
}(\mathcal{G},\theta)\longrightarrow C^{\ast}(\mathcal{G}_{F}\,,\theta
_{F})\longrightarrow0
\]
where $\theta_{U}$ (resp. $\theta_{F}$) is the 2-cocycle of $\mathcal{G}_{U}$
(resp. $\mathcal{G}_{F}$) obtained by restriction of $\theta$.
\end{prop}

This proposition was originally proven in \cite{Renault80} for $r$-discrete
groupoids and then mentioned in a more general setting in \cite{HilSkan87}. In
\cite{Ram98}, Proposition (\ref{exact}) was proven for $\theta=1$, and the
proof carries along immediately to the twisted case by using the full force of
\cite[Theorem 4.1]{Renault87} (see \cite[Part 1]{Latremoliere04}). The proof
of the upper semicontinuity in Theorem (\ref{cont}) follows from Proposition
(\ref{exact}) by observing that $p^{-1}(\{t\})\cap\mathcal{G}^{(0)}$ is a
saturated closed subset of $\mathcal{G}^{(0)}$ for all $t\in T$. The lower
semicontinuity relies on explicitly writing a family of well chosen $\theta
$-representations of $\mathcal{G}$, as detailed in \cite{Latremoliere04} and
\cite{Ram98} (for $\theta=1$).

Now, let us introduce the 2-cocycle $\gamma$ of $\Gamma$ by setting:%
\[
\gamma\left(  (n,\sigma,\chi),(n,\sigma,\chi^{\prime})\right)  =\sigma
(\chi,\chi^{\prime})\text{,}%
\]
for all $\left(  (n,\sigma,\chi),(n,\sigma,\chi^{\prime})\right)  \in
\Gamma^{(2)}$. It is straightforward to check that $\gamma$ is a 2-cocycle on
$\Gamma$. It remains to prove that $\gamma$ is continuous. Note that $\gamma$
lifts naturally to an (algebraic) 2-cocycle of $\Omega$, and $\gamma$ is
continuous on $\Gamma$ if and only if its lift to $\Omega$ is continuous since
the canonical surjection is open. Let us simply denote the lift of $\gamma$ to
$\Omega$ by $\gamma$ again. Since $\Omega$ is metrizable, it is enough to
prove that $\gamma$ is sequentially continuous. Let $((n_{j},\beta_{j}%
,\chi_{j}),(n_{j},\beta_{j},\chi_{j}^{\prime}))_{j\in\mathbb{N}}$ be a
sequence in $\Omega^{(2)}$ converging to $((n,\beta,\chi),(n,\beta
,\chi^{\prime}))$ in $\Omega^{(2)}$. Since $\widehat{G}$ is discrete, there
exists $N\in\mathbb{N}$ such that for all $j\geq N$ we have $(\chi_{j}%
,\chi_{,j}^{\prime})=(\chi,\chi^{\prime}).$ On the other hand, by definition
of the topology on $A\left(  \widehat{G}\right)  $, we then have
$\lim_{j\rightarrow\infty}\beta_{j}(\chi,\chi^{\prime})=\beta(\chi
,\chi^{\prime})$. Hence $\gamma((n_{j},\beta_{j},\chi_{j}),(n_{j},\beta
_{j},\chi_{j}^{\prime}))_{j\in\mathbb{N}}$ converges to $\gamma((n,\beta
,\chi),(n,\beta,\chi^{\prime}))$ and $\gamma$ is continuous.

\bigskip We shall denote by $((A_{t}:t\in T),B)$ the continuous field of
C*-algebras with parameter space $T$, where the continuous structure of the
family $(A_{t})_{t\in T}$ of C*-algebras is given by the C*-algebra $B$ of
continuous sections of $(A_{t})_{t\in T}$, as defined in \cite[Definition
10.3.1, p. 194]{Dixmier}. Since $\Gamma$ is amenable and separable by Lemma
(\ref{main_groupoid}), and we have proven that $\gamma$ is a continuous
2-cocycle of $\Gamma$, we have by Corollary (\ref{cont0}):

\begin{cor}
\label{cont_corr}$\left(  \left(  C^{\ast}\left(  \widehat{H_{n}}%
,\sigma\right)  :(n,\sigma)\in\Gamma^{(0)}\right)  ,C^{\ast}(\Gamma
,\gamma)\right)  $ is a continuous field of C*--algebras.
\end{cor}

Also, in exactly the same fashion, we can apply Corollary (\ref{cont0}) to the
bundle of transformation groupoids $\widetilde{\Gamma}$:

\begin{cor}
$\left(  \left(  C^{\ast}\left(  \widehat{H_{n}}\times_{\alpha_{n}}%
X_{n}\right)  :n\in\overline{\mathbb{N}}\right)  ,C^{\ast}\left(
\widetilde{\Gamma}\right)  \right)  $ is a continuous field of C*-algebras.
\end{cor}

\section{Quantum Gromov-Hausdorff convergence}

\subsection{Finite dimensional Approximations}

We denote by $\lambda_{n}$ the Haar probability measure on $H_{n}$ for all
$n\in\overline{\mathbb{N}}$, and we denote by $\left\Vert .\right\Vert
_{L^{1}}$ the norm of $L^{1}(H_{n},\lambda_{n})$. As a first step toward
Theorem (\ref{MAIN}), for each $n\in\overline{\mathbb{N}}$, we shall prove
that $C^{\ast}(\widehat{H_{n}},\sigma_{n})$ is the limit for
$\operatorname*{dist}_{q}$ of finite dimensional order-unit subspaces of
$C^{\ast}(\widehat{H_{n}},\sigma_{n})^{\operatorname*{sa}}$. These order-unit
spaces will be obtained, as in \cite{Rieffel00}, as the images of
multiplication operators defined on $C_{c}(\widehat{H_{n}})$ and extended to
$C^{\ast}(\widehat{H_{n}},\sigma_{n})$. In \cite{Rieffel00}, these operators
were constructed from an approximate unit for $L^{1}(\mathbb{T}^{d})$, known
as the Fejer kernels. Yet, this construction requires that $G$ is a Lie group.
For our purpose, we will use a slight modification of Rieffel's construction
which works for arbitrary compact groups but gives us a weaker approximation property:

\begin{lem}
\label{fejer}Let $f\in C(G)$ such that $f(e)=0.$ Let $\varepsilon>0$. There
exists a finite linear combination $\varphi\in C(G)$ of characters of $G$ such
that: $\varphi\geq0$, $\int_{G}\varphi d\lambda_{\infty}=1$ and $\int
_{G}\varphi(g)\left\vert f\right\vert (g)d\lambda_{\infty}(g)\leq\varepsilon$.
In this Lemma, $G$ can be any compact group.
\end{lem}

\begin{pf}
Let $\operatorname*{spec}(G)$ be the set of irreducible strongly continuous
representations of $G$, and let $\widehat{G}=\{\chi_{\pi}:\pi\in
\operatorname*{spec}(G)\}$ where $\chi_{\pi}$ is the character of the
representation $\pi$. Note that when $G$ is Abelian, $\widehat{G}%
=\operatorname*{spec}(G)$ is the dual of $G$. Yet for Theorem\ (\ref{Collapse}%
) it is useful to prove this lemma in full generality.

Let $U=\{g\in G:\left\vert f(g)\right\vert \leq\frac{1}{2}\varepsilon\}$,
which is an open neighborhood of $e$ by continuity of $f$ and since $f(e)=0$.
The family of open sets $\{g\in G:\chi(g)<\frac{1}{2}\}$ for $\chi\in
\widehat{G}$ covers the compact subset $G\backslash U$, so there exists a
finite set $F=\{\eta_{0},\ldots,\eta_{n}\}\subseteq\widehat{G}$ such that for
all $g\in G$ there exists $\eta_{j}\in F$ such that $\eta_{j}(g)<\frac{1}{2}$.
We can assume without loss of generality that the trivial representation is in
$F$. For each $\eta_{j}\in F$ let $\pi_{j}$ be the irreducible representation
of $G$ associated to $\eta_{j}$. Let $\pi=\oplus_{j=0}^{n}\pi_{j}$, and let
$\widetilde{\pi}=\pi\otimes\overline{\pi}$, where $\overline{\pi}$ is the
contragradiant representation of $\pi$. Let $\chi_{n}$ be the character of the
representation $\widetilde{\pi}^{\otimes n}$ for all $n\in\mathbb{N}$ (in
general, the representation $\widetilde{\pi}^{\otimes n}$ is not irreducible
so $\chi_{n}\not \in \widehat{G}$), and set $\varphi_{n}=\left\Vert \chi
_{n}\right\Vert _{L^{1}}^{-1}\chi_{n}$. Just as in \cite[Theorem
8.2]{Rieffel00}, we now claim for all $n\in\mathbb{N}$ that $\varphi_{n}\geq0$
and $\int_{G}\varphi_{n}d\lambda_{\infty}=1$, and above all we have that
$\int_{G\backslash U}\varphi_{n}d\lambda_{\infty}\rightarrow0$ when
$n\rightarrow\infty$. Therefore, there exists $N$ such that for all $n\geq N$
we have $\int_{G\backslash U}\varphi_{n}d\lambda_{\infty}<\varepsilon
(2\sup_{G}\left\vert f\right\vert )^{-1}$ and thus:%
\[
\int_{G}\left\vert f\right\vert \varphi_{n}d\lambda_{\infty}=\int
_{U}\left\vert f\right\vert \varphi_{n}d\lambda_{\infty}+\int_{G\backslash
U}\left\vert f\right\vert \varphi_{n}d\lambda_{\infty}\leq\frac{\varepsilon
}{2}\int_{G}\varphi_{n}d\lambda_{\infty}+\frac{\varepsilon}{2}\leq
\varepsilon\text{.}%
\]
On the other hand, by construction, $\chi_{1}$ is the finite sum of the
characters $\eta_{j}$ and $\overline{\eta_{j}}$ for $j\in\{0,\ldots,n\}$, so
$\chi_{n}$ is itself a finite linear combination of elements in $\widehat{G}$
and so is $\varphi_{n}$ for all $n\in\mathbb{N}$. Hence $\varphi=\varphi_{N}$
satisfies the conclusion of our Lemma. \qed \end{pf}

We wish to apply the following lemma to any of the groups $(H_{n}%
)_{n\in\overline{\mathbb{N}}}$, so we ease notations by working with a general
compact Abelian group $\Sigma$. We denote by $\beta$ a skew-bicharacter of
$\widehat{\Sigma}$. We denote again by $\alpha_{\Sigma}$ the dual action of
$\Sigma$ on $C^{\ast}\left(  \widehat{\Sigma},\beta\right)  $ and by
$\lambda_{\Sigma}$ the Haar probability measure on $\Sigma$. Given any
function $\varphi\in L^{1}(\Sigma)$, we define on $C_{c}(\widehat{\Sigma})$
the operator:%
\[
\alpha_{\Sigma}^{\varphi}:f\in C_{c}(\widehat{\Sigma})\longmapsto\int_{\Sigma
}\varphi(g)\alpha_{\Sigma}^{g^{-1}}(f)d\lambda_{\Sigma}(g)\in C_{c}%
(\widehat{\Sigma})\text{.}%
\]
Since $\Sigma$ acts by *-automorphisms, hence by isometries, on $C^{\ast
}(\widehat{\Sigma},\beta)$, it is immediate that for all $f\in C_{c}%
(\widehat{\Sigma})$ we have $\left\Vert \alpha_{\Sigma}^{\varphi
}(f)\right\Vert _{C^{\ast}\left(  \widehat{\Sigma},\beta\right)  }%
\leq\left\Vert \varphi\right\Vert _{L^{1}}\left\Vert f\right\Vert _{C^{\ast
}\left(  \widehat{\Sigma},\beta\right)  }$ where we denote $\int\left\vert
\varphi\right\vert d\lambda_{\Sigma}$ by $\left\Vert \varphi\right\Vert
_{L^{1}}$. By abuse of notation, we also denote by $\alpha_{\Sigma}^{\varphi}$
the unique continuous extension of $\alpha_{\Sigma}^{\varphi}$ to $C^{\ast
}(\widehat{\Sigma},\beta)$. We summarize in the following lemma some easy
properties of the operator $\alpha_{\Sigma}^{\varphi}$ when $\varphi$ is a
linear combination of characters of $\Sigma$, in view of an application of
Lemma (\ref{alpha_phi.1}) and Lemma\ (\ref{fejer}) together:

\begin{lem}
\label{alpha_phi.1}Let $\varphi$ be a linear combination of characters of
$\Sigma$. The range of $\alpha_{\Sigma}^{\varphi}$ is the finite dimensional
space:%
\[
R(\varphi)=\left\{  f\in C_{c}(\widehat{\Sigma}):\operatorname*{supp}%
f\subseteq\operatorname*{supp}(\widehat{\varphi})\right\}
\]
where $\widehat{\varphi}$ is the Fourier transform of $\varphi.$In particular,
the dimension of $R(\varphi)$ is the cardinal of $\operatorname*{supp}%
(\widehat{\varphi})$. Moreover, if $\varphi$ is real-valued then the operator
$\alpha_{\Sigma}^{\varphi}$ is self-adjoint and therefore $\alpha_{\Sigma
}^{\varphi}\left(  C^{\ast}\left(  \widehat{\Sigma},\beta\right)
^{\operatorname*{sa}}\right)  $ is the finite dimensional order-unit subspace
$V(\varphi)=R(\varphi)^{\operatorname*{sa}}$ of $C^{\ast}(\widehat{\Sigma
},\beta)^{\operatorname*{sa}}$.
\end{lem}

\begin{pf}
By definition, $\varphi\in L^{1}(\Sigma)$ and thus the operator $\alpha
_{\Sigma}^{\varphi}$ is well defined and of norm $\int_{\Sigma}\left\vert
\varphi\right\vert d\lambda_{\Sigma}$. Let $f\in C_{c}\left(  \widehat{\Sigma
},\beta\right)  \subseteq C^{\ast}\left(  \widehat{\Sigma},\beta\right)  $.
Let us compute $\alpha_{\Sigma}^{\chi}\left(  C_{c}\left(  \widehat{\Sigma
}\right)  \right)  $ for a character $\chi\in\widehat{\Sigma}$. Let
$\chi^{\prime}\in\widehat{\Sigma}$. Denote by $\delta_{\chi}$ the Dirac
measure at $\chi$. We observe that:%
\[
\alpha_{\Sigma}^{\chi}(f)(\chi^{\prime})=\int_{\Sigma}\chi(g)\chi^{\prime
}(g^{-1})f(\chi^{\prime})d\lambda_{\Sigma}(g)=f(\chi^{\prime})\int_{\Sigma
}\chi\overline{\chi^{\prime}}d\lambda_{\Sigma}\text{,}%
\]
so $\alpha_{\Sigma}^{\chi}(f)(\chi^{\prime})=f(\chi^{\prime})\delta_{\chi
}(\chi^{\prime})$, hence $\alpha_{\Sigma}^{\chi}\left(  C_{c}(\widehat{\Sigma
})\right)  $ is the one dimensional subspace $R(\chi)$ of $C^{\ast}\left(
\widehat{\Sigma},\beta\right)  $.

By linearity of $\varphi\mapsto\alpha_{\Sigma}^{\varphi}$, we deduce that
$R(\varphi)=\alpha_{\Sigma}^{\varphi}\left(  C_{c}\left(  \widehat{\Sigma
}\right)  \right)  $ for any linear combination $\varphi$ of characters of
$\Sigma$. Yet, since $\alpha_{\Sigma}^{\varphi}$ is continuous and
$C_{c}\left(  \widehat{\Sigma}\right)  $ is dense in $C^{\ast}\left(
\widehat{\Sigma},\beta\right)  $, the range of $\alpha_{\Sigma}^{\varphi}$ is
included in the closure of $R(\varphi)$, which is $R(\varphi)$ since
$R(\varphi)$ is finite dimensional.

Assume now that $\varphi$ is real valued. Since $\beta$ is a skew bicharacter
of $\widehat{\Sigma}$, so that $\beta(\chi,\chi)=1$ for all $\chi\in\Sigma$,
we check that:%
\begin{align}
\alpha_{\Sigma}^{\varphi}(f^{\ast})(\chi)  &  =\int_{\Sigma}\varphi
(g)\chi(g^{-1})\beta(\chi,\chi)\overline{f(\chi^{-1})}d\lambda_{\Sigma
}(g)\label{alpha_phi.1.1}\\
&  =\int_{\Sigma}\overline{\varphi(g)\chi(g)f(\chi^{-1})}d\lambda_{\Sigma
}(g)=\overline{\alpha_{\Sigma}^{\varphi}(f)(\chi^{-1})}=\left(  \alpha
_{\Sigma}^{\varphi}(f)\right)  ^{\ast}(\chi)\text{.}%
\end{align}
By continuity of $\alpha_{\Sigma}^{\varphi}$ and density of
$C_{c}\left(  \widehat{\Sigma},\beta\right)  $ in $C^{\ast}\left(  \widehat{\Sigma
},\beta\right)  $, Equality (\ref{alpha_phi.1.1}) extends to $C^{\ast}\left(  \widehat{\Sigma
},\beta\right)  $, so $\alpha_{\Sigma}^{\varphi}$ is selfadjoint on $C^{\ast}\left(  \widehat{\Sigma},\beta\right)  $. The image
of $C^{\ast}\left(  \widehat{\Sigma},\sigma\right)  ^{\operatorname*{sa}}$ is
deduced immediately from the self-adjointness of $\alpha_{\Sigma}^{\varphi}$. \qed \end{pf}

\bigskip We now return to our specific setting. Lemma\ (\ref{fejer}) gives us,
for any $\varepsilon>0$, a finite linear combination of characters
$\varphi_{\varepsilon}$ of $G$ which we can use to define the operator
$\alpha_{G}^{\varphi_{\varepsilon}}$ which has the properties given in Lemma
(\ref{alpha_phi.1}). Yet for all $n\in\mathbb{N}$, we can also define the
operators $\alpha_{H_{n}}^{\varphi_{\varepsilon}}$ on $C^{\ast}\left(
\widehat{H_{n}},\sigma_{n}\right)  $, by using the restriction of
$\varphi_{\varepsilon}$ to $H_{n}$. Since $H_{n}$ is a closed subgroup of $G$,
the restriction of $\varphi_{\varepsilon}$ to $H_{n}$ is a linear combination
of characters of $H_{n}$. The operators $\alpha_{H_{n}}^{\varphi_{\varepsilon
}}$ satisfy the properties of Lemma (\ref{alpha_phi.1}), and thus their ranges
are finite dimensional. The following two lemmas are the key step which will
allow us to identify the ranges of all the operators $\alpha_{H_{n}}%
^{\varphi_{\varepsilon}}$ for $n\in\overline{\mathbb{N}}$ and $n$ large enough.

\begin{lem}
\label{finite1 copy(1)}Let $F$ be a finite subset of $\widehat{G}%
\backslash\{1\}$. There exists $N\in\mathbb{N}$ such that for all $n\geq N$ we
have $J_{n}\cap F=\varnothing$.
\end{lem}

\begin{pf}
Assume that for all $N\in\mathbb{N}$ there exists $n\geq N$ such that
$J_{n}\cap F\neq\varnothing$. This implies $F$ intersects countably many of
the groups $J_{n}$, and since $F$ is finite, we deduce that there exists an
element $\chi\in F$ such that, for some strictly increasing function
$f:\mathbb{N\rightarrow N}$ we have: $\chi\in\bigcap_{n\in\mathbb{N}}J_{f(n)}%
$. Yet, this implies that the character $\chi$ of $G$ is $1$ on $\bigcup
_{n\in\mathbb{N}}H_{f(n)}$. By Lemma (\ref{dense}), $\bigcup_{n\in\mathbb{N}%
}H_{f(n)}$ is dense in $G$ since $\lim_{\infty}f=\infty$, so $\chi=1$. This
contradicts the hypothesis that $F$ is a subset of $\widehat{G}\backslash
\{1\}$. We have proven our lemma by contradiction. \qed \end{pf}

\begin{rem}
If $\widehat{G}$ admits is a proper space for some metric (for instance if
$\widehat{G}$ is finitely generated), Lemma (\ref{finite1 copy(1)}) proves
that all balls in $J_{n}$ around $1\in\widehat{G}$, which are then finite, are
in fact all equal to $\{1\}$ for $n$ large enough. This proves that
$(J_{n})_{n\in\mathbb{N}}$ converges to $\{1\}$ for the Gromov-Hausdorff
topology on locally compact spaces \cite[Definition 3.14]{Gromov}, which is
the expected dual statement to the convergence of $(H_{n})_{n\in\mathbb{N}}$.
\end{rem}

\begin{cor}
\label{finiteembedd}Let $F$ be a finite subset of $\widehat{G}$. There exists
$N\in\mathbb{N}$ such that for all $n\geq N$, the canonical surjection
$q_{n}:\widehat{G}\rightarrow\widehat{H_{n}}$ is injective on $F$.
\end{cor}

\begin{pf}
Let $F^{\prime}=\{\chi^{-1}\chi^{\prime}\in G:\chi,\chi^{\prime}\in
F\}\backslash\{1\}$. Then $F^{\prime}$ is finite since $F$ is. By Lemma
(\ref{finite1 copy(1)}), there exists $N\in\mathbb{N}$ such that for all
$n\geq N$ we have $J_{n}\cap F^{\prime}=\varnothing$. Let $n\geq N$. Now, if
$q_{n}(\chi)=q_{n}(\chi^{\prime})$ then $\chi^{-1}\chi^{\prime}\in J_{n}$, so
$\chi^{-1}\chi^{\prime}\not \in F^{\prime}$, so if $\chi,\chi^{\prime}\in F$
then $\chi^{-1}\chi^{\prime}=1$, and our lemma is proved. \qed \end{pf}

The following lemma will also be useful to deal with matters of normalization
of some operators. We include our proof for the reader's convenience. We 
wish to thank the referee for providing the following references for this 
result: \cite[Theorem in Annex]{Glimm62}, 
\cite[Ch. 8, sec. 5,  n$^{\text{0}}$ 6]{BourbakiIntegration}.

\begin{lem}
\label{riemann}Let $f\in C(G),$ where $C(G)\ $is the C*-algebra of continuous
function on $G$. Then:%
\[
\lim_{n\rightarrow\infty}\int_{H_{n}}fd\lambda_{n}=\int_{G}fd\lambda_{\infty
}\text{.}%
\]

\end{lem}

\begin{pf}
We shall prove that every subsequence of $(\lambda_{n})_{n\in\mathbb{N}}$
admits a subsequence weakly converging to $\lambda_{\infty}$, which will then
prove our Lemma. Let $(\lambda_{m(n)})_{n\in\mathbb{N}}$ be a subsequence
$(\lambda_{n})_{n\in\mathbb{N}}$. Since $G$ is compact, $(\lambda_{m(n)})_{n\in\mathbb{N}}$ is uniformly tight \cite[Definition p.
293, Theorem 11.5.4]{Dudley}\ and so there exists a weakly convergent
subsequence $(\lambda_{k(n)})_{n\in\mathbb{N}}$ of $(\lambda_{m(n)}%
)_{n\in\mathbb{N}}$. Let $I:f\in C(G)\mapsto\lim_{n\rightarrow\infty}%
\int_{H_{k(n)}}fd\lambda_{k(n)}$. The map $I$ is obviously a continuous linear
functional on $C(G)$, and $I(1)=1$. To prove our lemma, it is enough to show
that $I$ is also invariant by translation, as the conclusion will then follows
from the uniqueness of the Haar probability measure on $G$. Let $f\in C(G)$
and let $g\in G$. Let $\varepsilon>0$. Since $f$ is continuous on the compact
$G$, it is uniformly continuous so there exists $\delta>0$ such that for all
$g^{\prime},g^{\prime\prime}\in G$ such that $l(g^{\prime-1}g^{\prime\prime
})<\delta$ we have $\left\vert f(g^{\prime\prime})-f(g^{\prime})\right\vert
<\varepsilon$. Since $\left(  H_{k(n)}\right)  _{n\in\mathbb{N}}$ converges to
$G$ for $\mathfrak{H}(l)$, there exists $N\in\mathbb{N}$ such that for all
$n\geq N$ there exists $g_{n}\in H_{k(n)}$ such that $l(g_{n}^{-1}g)<\delta$.
We now have for all $n\geq N$ and all $x\in G$:%
\begin{align*}
&  \left\vert \int_{H_{k(n)}}\left(  f(xg)-f(x)\right)  d\lambda
_{k(n)}(x)\right\vert \leq\\
&  \leq\int_{H_{k(n)}}\left\vert f(xg_{n})-f(xg)\right\vert d\lambda
_{k(n)}(x)+\left\vert \int_{H_{k(n)}}(f(xg_{n})-f(x))d\lambda_{k(n)}%
(x)\right\vert \\
&  \leq\varepsilon+\left\vert \int_{H_{k(n)}}(f(xg_{n})-f(x))d\lambda
_{k(n)}(x)\right\vert
\end{align*}
using uniform continuity of $f$ and $l\left(  \left(  xg_{n}\right)
^{-1}xg\right)  =l(g_{n}^{-1}g)<\delta$. 
Now, since $\lambda_{k(n)}$ is translation invariant
on $H_{k(n)}$, we have $\int_{H_{k(n)}}(f(xg_{n})-f(x))d\lambda_{k(n)}(x)=0$. Taking the limit, we have proven that, if $f_{g}:x\mapsto
f(xg)$, then for all $\varepsilon>0$ and all $g\in G$ we have $\left\vert
I(f_{g})-I(f)\right\vert <\varepsilon$, so $I(f_{g})=I(f)$ so our Lemma is proven. \qed \end{pf}

We are now in a position to prove that for the quantum Gromov-Hausdorff
distance, for any $n\in\overline{\mathbb{N}}$, the quantum metric space
$C^{\ast}(\widehat{H_{n}},\sigma)$ is the limit of spaces of the form
$V(\varphi)$, with the notation of Lemma (\ref{alpha_phi.1}). We will use the
following result from Rieffel, which we quote without proof for convenience:

\begin{lem}
\label{rieffel}Let $(A,L)$ be a compact quantum metric space and $B\subseteq
A$ an order unit subspace of $A$. Assume there exists a map $P:A\rightarrow B$
and $\delta\in]0,+\infty)$ such that for all $a\in A$ we have $L(P(a))\leq
L(a)$ and $\left\Vert a-P(a)\right\Vert \leq\delta L(a)$. Then
$\operatorname*{dist}_{q}((A,L),(B,L))\leq\delta$.
\end{lem}

In our context, the map $P$ is of the form $\alpha_{\Sigma}^{\varphi}$ for a
well-chosen $\varphi$. Now, we can prove the following proposition:

\begin{prop}
\label{finitedimapprox}Let $\varepsilon>0$. There exists $N_{\varepsilon}%
\in\mathbb{N}$ and a finite dimensional space $V_{\varepsilon}$ such that, for
all $n\in\overline{\mathbb{N}}$ such that $n\geq N_{\varepsilon}$, there
exists a linear isomorphism $\theta_{n}:V_{\varepsilon}\longrightarrow
C^{\ast}\left(  \widehat{H_{n}},\sigma_{n}\right)  ^{\operatorname*{sa}}$ such
that $\theta_{n}(V)$ is an order-unit subspace of $C^{\ast}\left(
\widehat{H_{n}},\sigma_{n}\right)  ^{\operatorname*{sa}}$ and:%
\[
\operatorname*{dist}\nolimits_{q}\left(  \left(  \theta_{n}(V_{\varepsilon
}),L_{n}\right)  ,\left(  C^{\ast}\left(  \widehat{H_{n}},\sigma_{n}\right)
,L_{n}\right)  \right)  \leq\frac{1}{3}\varepsilon\text{.}%
\]
Explicitly, we can choose $\theta_{\infty}=\operatorname*{Id}$ (so
$V_{\varepsilon}\subseteq C^{\ast}(\widehat{G},\sigma_{\infty})$) and for each
$n\in\mathbb{N}$, we can set $\theta_{n}$ to be the restriction to
$V_{\varepsilon}$ of the integrated quotient map $\Theta_{n}$ defined for all
$f\in C_{c}(\widehat{G})$ to be the map $\Theta_{n}(f):q_{n}(h)\in
\widehat{H_{n}}\mapsto\sum_{j\in J_{n}}f(hj)$.
\end{prop}

\begin{pf}
By Lemma\ (\ref{fejer}) applied with $f$ as the length function $l$, there
exists a finite linear combination $\varphi_{\varepsilon}$ of characters of
$G$ such that:%
\[
\int_{G}\varphi_{\varepsilon}(g)l(g)d\lambda_{\infty}(g)\leq\frac{\varepsilon
}{3(1+\varepsilon)}\text{.}%
\]

Let $B=\operatorname*{supp}(\widehat{\varphi_{\varepsilon}})\subseteq
\widehat{G}$. Let $V_{\varepsilon}=V(\varphi_{\varepsilon})$ as in Lemma
(\ref{fejer}), which is the finite dimensional order-unit space $\left\{  f\in
C_{c}(\widehat{G})^{\operatorname*{sa}}:\operatorname*{supp}(f)\subseteq
B\right\}  $, namely the self-adjoint part of the range of $\alpha
_{G}^{\varphi_{\varepsilon}}$ (for the involution in $C^{\ast}(\widehat
{G},\sigma_{\infty})$). Now, for any $n\in\mathbb{N}$, the restriction of
$\varphi_{\varepsilon}$ to $H_{n}$ is a linear combination of characters of
$\widehat{H_{n}}$ and therefore, by Lemma (\ref{fejer}), the self-adjoint part
of the range of $\alpha_{H_{n}}^{\varphi_{\varepsilon}}$ is the order-unit
subspace $V_{\varepsilon}^{n}=\left\{  f\in C_{c}\left(  \widehat{H_{n}%
}\right)  ^{\operatorname*{sa}}:\operatorname*{supp}(f)\subseteq
q_{n}(B)\right\}  $ of $C^{\ast}\left(  \widehat{H_{n}},\sigma_{n}\right)  $.

We set $\theta_{\infty}$ to be the identity map on $V_{\varepsilon}$. By
construction in Lemma (\ref{fejer}), $B$ is a finite set, so by Corollary
(\ref{finiteembedd}) there exists $N_{1}\in\mathbb{N}$ such that, for all
$n\in\mathbb{N}$ such that $n\geq N_{1}$ the canonical surjection
$q_{n}:\widehat{G}\rightarrow\widehat{H_{n}}$ is injective on $B$. Therefore,
the order-unit subspace $V_{\varepsilon}^{n}$ is of the same dimension as
$V_{\varepsilon}$ for $n\geq N_{1}$. This proves already the existence of the
linear isomorphism $\theta_{n}$. It is however useful for our argument to
choose a specific map $\theta_{n}$. Let $\theta_{n}$ be the restriction to
$V_{\varepsilon}$ of the map $\Theta_{n}$. It is straightforward that
$\Theta_{n}$ is a well-defined linear surjection, and that
$\operatorname*{supp}(\Theta_{n}(f))\subseteq q_{n}(\operatorname*{supp}(f)).$
In particular, $\Theta_{n}(V_{\varepsilon})=V_{\varepsilon}^{n}$ and, since
both $V_{\varepsilon}$ and $V_{\varepsilon}^{n}$ have the same finite
dimension, $\Theta_{n}$ is indeed a linear isomorphism from $V_{\varepsilon}$
onto $V_{\varepsilon}^{n}$.

By Lemma (\ref{riemann}), we have $\int_{H_{n}}\varphi d\lambda_{n}%
\underset{n\rightarrow\infty}{\longrightarrow}\int_{G}\varphi d\lambda=1$, so
there exists $N_{2}\in\mathbb{N}$ such that for all $n\geq N_{2}$ we have
$\int_{H_{n}}\varphi_{\varepsilon}d\lambda_{n}\geq\left(  1+\varepsilon
\right)  ^{-1}>0$. Set $N_{\varepsilon}=\max\{N_{1},N_{2}\}$.For all
$n\in\overline{\mathbb{N}}$ such that $n\geq N_{\varepsilon}$, we set:%
\[
c_{n}=\left(  \int_{H_{n}}\varphi_{\varepsilon}(h)d\lambda_{n}(h)\right)
^{-1}=\left\Vert \alpha_{H_{n}}^{\varphi}\right\Vert ^{-1}%
\]
and we note that $c_{n}\leq1+\varepsilon$. By lemma (\ref{alpha_phi.1})
applied with $\Sigma=H_{n}$, the operator $c_{n}\alpha_{H_{n}}^{\varphi
_{\varepsilon}}$ is a linear map of $C^{\ast}\left(  \widehat{H_{n}}%
,\sigma_{n}\right)  $ which maps $C^{\ast}\left(  \widehat{H_{n}},\sigma
_{n}\right)  ^{\operatorname*{sa}}$ onto the finite dimensional order-unit
subspace $V_{\varepsilon}^{n}(\varphi_{\varepsilon})$ (since $c_{n}%
\in\mathbb{R}$). We wish to apply Lemma (\ref{rieffel}) with $C^{\ast}\left(
\widehat{H_{n}},\sigma_{n}\right)  $ for $A$, the subspace $V_{\varepsilon
}^{n}(\varphi_{\varepsilon})$ for $B$ and the operator $c_{n}\alpha_{H_{n}%
}^{\varphi_{\varepsilon}}$ for the function $P.$

Let $a\in C^{\ast}\left(  \widehat{H_{n}},\sigma_{n}\right)  $. \ Since
$c_{n}\int_{H_{n}}\varphi_{\varepsilon}d\lambda_{n}=1$, we have:%
\begin{align}
\left\Vert a-c_{n}\alpha_{H_{n}}^{\varphi_{\varepsilon}}(a)\right\Vert _{n}
&  =\left\Vert \int_{H_{n}}c_{n}\varphi_{\varepsilon}(g)(a-\alpha_{H_{n}%
}^{g^{-1}}(a))d\lambda_{n}(g)\right\Vert _{n}\nonumber\\
&  \leq\int_{H_{n}}c_{n}\varphi_{\varepsilon}(g)\left\Vert a-\alpha_{H_{n}%
}^{g^{-1}}(a)\right\Vert _{n}d\lambda_{n}(g)\nonumber\\
&  \leq c_{n}\left(  \int_{H_{n}}\varphi_{\varepsilon}(g)l(g)d\lambda
_{n}(g)\right)  L_{n}(a)\leq\frac{\varepsilon c_{n}L_{n}(a)}{3(1+\varepsilon)}
\label{finitedimapprox.1}%
\end{align}

since by definition we have for all $g\in H_{n}$ that
$L_{n}(a)l(g)\geq\left\Vert a-\alpha_{H_{n}}^{g}(a)\right\Vert _{n}$. We also
observe that $\alpha_{H_{n}}^{g}$ and $\alpha_{H_{n}}^{\varphi_{\varepsilon}}$
commute for all $g\in H_{n}$, so 
\[ 
\left\Vert \alpha_{H_{n}}^{g}(c_{n}%
\alpha_{H_{n}}^{\varphi_{\varepsilon}}(a))-c_{n}\alpha_{H_{n}}^{\varphi
_{\varepsilon}}(a)\right\Vert _{n}\leq\left\Vert \alpha_{H_{n}}^{g}%
(a)-a\right\Vert _{n}\] since $c_{n}\alpha_{H_{n}}^{\varphi_{\varepsilon}}$ is
of norm 1. Hence:
\begin{equation}
L_{n}(c_{n}\alpha_{H_{n}}^{\varphi_{\varepsilon}}(a))\leq L_{n}%
(a).\label{finitedimapprox.2}%
\end{equation}

Hence, with $\delta=\frac{1}{3}c_{n}\varepsilon$, using
(\ref{finitedimapprox.1}) and (\ref{finitedimapprox.2}) in Lemma
(\ref{rieffel}), we get:%
\[
\operatorname*{dist}\nolimits_{q}\left(  \left(  \theta_{n}(V_{\varepsilon
}),L_{n}\right)  ,\left(  C^{\ast}\left(  \widehat{H_{n}},\sigma_{n}\right)
,L_{n}\right)  \right)  \leq\frac{\varepsilon}{3(\varepsilon+1)}c_{n}\leq
\frac{1}{3}\varepsilon\text{,}%
\]
since $n\geq N_{\varepsilon}\geq N_{2}$ we have $c_{n}\leq1+\varepsilon$. Our
proof is concluded. \qed \end{pf}

In what follows, we choose $\theta_{n}$ to be the restriction of the map
$\Theta_{n}$ to $V_{\varepsilon}$ for all $n\in\overline{\mathbb{N}}$.

\subsection{Continuous Field of Lip-norms}

We fix $\varepsilon>0$. By Proposition (\ref{finitedimapprox}), there exists
$N_{\varepsilon}\in\mathbb{N}$ and $\varphi_{\varepsilon}\in C(G)$ such that,
if $V_{\varepsilon}=V(\varphi_{\varepsilon})$, the spaces $\theta
_{n}(V(\varphi_{\varepsilon}))$ approximates $C^{\ast}(\widehat{H_{n}}%
,\sigma_{n})$ within $\varepsilon/3$ for $\operatorname*{dist}_{q}$ for all
$N_{\varepsilon}\leq n\leq\infty$. We now have reduced our problem to proving
that the family of finite dimensional quantum metric spaces $(\theta
_{n}(V_{\varepsilon}),L_{n})_{n\in\{N_{\varepsilon},N_{\varepsilon
}+1,...,\infty\}}$ forms a continuous field of compact quantum metric spaces.
The method is similar to the one in \cite{Rieffel00}, yet it is adapted to the
fact we are dealing here with many groups rather than just one.

The first observation is a direct consequence of the first part of our paper:

\begin{lem}
\label{Vcontnorm}Let $a\in V_{\varepsilon}$. The function $n\in\overline
{\mathbb{N}}\mapsto\left\Vert \theta_{n}(a)\right\Vert _{n}$ is continuous (at infinity).
\end{lem}

\begin{pf}
Let $a\in V_{\varepsilon}$. We define $f\in C_{c}(\Gamma)$ as
follows. Let $(n,\sigma,\chi)\in\Gamma$. If $n<N_{\varepsilon}$ we set
$f(n,\sigma,\chi)=0$. Assume now that $n\geq N_{\varepsilon}$. We denote by
$q_{n}$ the canonical surjection $\widehat{G}\rightarrow\widehat{H_{n}}$. If
$\chi\not \in q_{n}(\operatorname*{supp}(\widehat{\varphi_{\varepsilon}}))$
then we set $f(n,\sigma,\chi)=0$. If $\chi\in q_{n}(\operatorname*{supp}%
(\widehat{\varphi_{\varepsilon}}))$ then there exists a unique $\chi^{\prime
}\in\widehat{G}$ such that $q_{n}(\chi^{\prime})=\chi$ by corollary
(\ref{finiteembedd}). We set $f(n,\sigma,\chi)=a(\chi^{\prime})$. The function
$f$ is continuous by construction.

The support of $f$ is contained in the set $\Gamma^{(0)}\times
\operatorname*{supp}(\widehat{\varphi_{\varepsilon}})$ which is compact as
$\operatorname*{supp}(\widehat{\varphi_{\varepsilon}})$ is finite and
$\Gamma^{(0)}$ is compact by construction. Therefore, $f$ defines an element
in $C^{\ast}(\Gamma,\gamma)$. By Corollary (\ref{cont_corr}), $n\longmapsto
\left\Vert f^{(n,\sigma_{n})}\right\Vert _{n}$ is continuous on $\overline
{\mathbb{N}}$ (since we assumed that $\lim_{n\rightarrow\infty}\sigma
_{n}=\sigma_{\infty}$).

For $n\geq N_{\varepsilon}$, by choosing in Proposition (\ref{finitedimapprox}%
) the map $\theta_{n}$ to be the restriction of the map $\Theta_{n}$, it is
immediate that $f(n,\sigma,\chi)=\theta_{n}(a)(\chi)$ for $n\geq
N_{\varepsilon}$. We therefore deduce that $f^{(n,\sigma_{n})}=\theta_{n}(a)$
and our lemma is proven. \qed \end{pf}

Now, we turn to the continuity of the Lip-norms. For all $a\in V_{\varepsilon
}$ and $g\in G\backslash\{e\}$, we set $\Delta_{g}(a)=l(g)^{-1}\left(
\alpha_{G}^{g}(a)-a\right)  $.

Observe that $\Delta_{g}(a)\in V_{\varepsilon}$ as $\alpha_{G}$ restricts to
an action on $V_{\varepsilon}$, so $L_{n}(\theta_{n}(a))=\sup\left\{
\left\Vert \theta_{n}(\Delta_{g}(a))\right\Vert _{n}:g\in H_{n}\backslash
\{e\}\right\}  $ since $\theta_{n}(a)\in C_{c}(\widehat{H_{n}})$, and so
$\theta_{n}(\alpha_{G}^{g}(a))=\alpha_{H_{n}}^{g}(\theta_{n}(a))$ for any
$g\in H_{n}$ for our choice of $\theta_{n}=\Theta_{n}$ on $V_{\varepsilon}$.

For $a\in V_{\varepsilon}$, the function $n\in\overline{\mathbb{N}}%
\mapsto\left\Vert \theta_{n}(a)\right\Vert _{n}$ is continuous on the compact
$\overline{\mathbb{N}}\backslash\{0,1,\ldots,N_{\varepsilon}-1\}$ by Lemma
(\ref{Vcontnorm}), so it is bounded and so we can define $\left\Vert
a\right\Vert _{\ast}=\sup_{n\in\overline{\mathbb{N}}\backslash\{0,1,\ldots
,N_{\varepsilon}-1\}}\left\Vert \theta_{n}(a)\right\Vert _{n}$. It is
straightforward to check that $\left\Vert .\right\Vert _{\ast}$ is a norm on
$V_{\varepsilon}$.

\begin{prop}
Let $a\in V_{\varepsilon}$. The function $n\in\overline{\mathbb{N}}\mapsto
L_{n}(\theta_{n}(a))$ is continuous (at infinity).
\end{prop}

\begin{pf}
Let $\varepsilon>0$ and $a\in V_{\varepsilon}$. There exists $g_{\varepsilon
}\in G\backslash\{e\}$ such that $L_{\infty}(a)-\varepsilon\leq\left\Vert
\Delta_{g_{\varepsilon}}(a)\right\Vert _{\infty}\leq L_{\infty}(a)$. Now, the
function $g\mapsto\Delta_{g}(a)$ is continuous at $g_{\varepsilon}\in G$ for
the norm $\left\Vert .\right\Vert _{\ast}$. Indeed, the dual action
$\alpha_{G}$ is continuous on $(V_{\varepsilon},\left\Vert .\right\Vert
_{\infty})$, and since $V_{\varepsilon}$ is finite dimensional, all the norms
on $V_{\varepsilon}$ are equivalent. Moreover, by assumption $l$ is continuous
and $l(g)=0$ only when $g=e$, so the function $g\in G\backslash\{e\}\mapsto
l(g)^{-1}$ is continuous as well.

Therefore, we can find $\delta>0$ such that, for all $g\in G$ such that
$l(g_{\varepsilon}g^{-1})\leq\delta$, we have:%
\[
\left\Vert \Delta_{g}(a)-\Delta_{g_{\varepsilon}}(a)\right\Vert _{\ast}%
\leq\varepsilon
\]
by continuity of $g\in G\backslash\{e\}\mapsto\Delta_{g}(a)$ for $\left\Vert
.\right\Vert _{\ast}$ at $g_{\varepsilon}$.

Also, by Lemma (\ref{Vcontnorm}), we can choose $N\in\mathbb{N}$ so that for
all $n\geq N$:
\[
\left\vert \left\Vert \theta_{n}(\Delta_{g_{\varepsilon}}(a))\right\Vert
_{n}-\left\Vert \Delta_{g_{\varepsilon}}(a)\right\Vert _{\infty}\right\vert
\leq\varepsilon.
\]
On the other hand, since $(H_{n})_{n\in\mathbb{N}}$ converges for the
Hausdorff distance $\mathfrak{H}(l)$ to $G$, there exists $N^{\prime}%
\in\mathbb{N}$ such that for all $n\geq N^{\prime}$, any point in $G$ is
within an $\delta$-ball whose center is in $H_{n}$. In particular, we can
find, for all $n\geq N^{\prime}$, an element $g_{n}\in H_{n}$ such that
$l(g_{\varepsilon}g_{n}^{-1})\leq\delta.$ We let $N^{\prime\prime}%
=\max(N,N^{\prime})$.

We can now compute, for all $n\geq N^{\prime\prime}$:
\begin{align*}
\left\vert \left\Vert \theta_{n}(\Delta_{g_{n}}(a))\right\Vert _{n}-\left\Vert
\Delta_{g_{\varepsilon}}(a)\right\Vert _{\infty}\right\vert  &  \leq\left\vert
\left\Vert \theta_{n}(\Delta_{g_{n}}(a))\right\Vert _{n}-\left\Vert \theta
_{n}(\Delta_{g_{\varepsilon}}(a))\right\Vert _{n}\right\vert \\
&  +\left\vert \left\Vert \theta_{n}(\Delta_{g_{\varepsilon}}(a))\right\Vert
_{n}-\left\Vert \Delta_{g_{\varepsilon}}(a)\right\Vert _{\infty}\right\vert \\
&  \leq\left\Vert \theta_{n}(\Delta_{g_{n}}(a)-\Delta_{g_{\varepsilon}%
}(a))\right\Vert _{n}+\varepsilon\\
&  \leq\left\Vert \Delta_{g_{n}}(a)-\Delta_{g_{\varepsilon}}(a)\right\Vert
_{\ast}+\varepsilon\leq2\varepsilon
\end{align*}
where we used that $\theta_{n}$ is a linear map of norm 1 as the restriction
of $\Theta_{n}$ by Proposition (\ref{finitedimapprox}). Now we can conclude
that, for $n\geq N^{\prime\prime}$:
\[
L_{n}(\theta_{n}(a))\geq\left\Vert \theta_{n}(\Delta_{g_{n}}(a))\right\Vert
_{n}\geq\left\Vert \Delta_{g_{\varepsilon}}(a)\right\Vert _{\infty
}-2\varepsilon\geq L_{\infty}(a)-3\varepsilon\text{.}%
\]
Hence $\lim\inf_{n\rightarrow\infty}L_{n}(\theta_{n}(a))\geq L_{\infty}(a)$.

On the other hand, it is obvious that $L_{n}(\theta_{n}(a))\leq L_{n}^{\prime
}(a)$ where $L_{n}^{\prime}(a)=\sup\{\left\Vert \theta_{n}(\Delta
_{g}(a))\right\Vert _{n}:g\in G\backslash\{e\}\}$. Yet, we are now in the
setting of \cite[Lemma 9.3]{Rieffel00}: we have a continuous field of norms on
a finite dimensional vector space, and we define the Lip-norms $L_{n}^{\prime
}$ as the supremum, for each of those norms, of the same sets (because now we
only deal with one group, $G$). We therefore omit this part and we refer the
reader to \cite[Lemma 9.3]{Rieffel00} or \cite{Latremoliere04} for the details. \qed \end{pf}

\bigskip In fact, we have proven a little more. Let $(H_{n}^{\prime}%
)_{n\in\mathbb{N}}$ be a sequence of closed subgroups of $G$ such that
$H_{n}\subseteq H_{n}^{\prime}$ for all $n\in\overline{\mathbb{N}}$. We can
define a sequence of Lip-norms \
\begin{equation}
L_{H_{n}^{\prime}}:a\in V_{\varepsilon}\longmapsto\sup\left\{  \left\Vert
\theta_{n}(\Delta_{g}(a))\right\Vert _{n}:g\in H_{n}^{\prime}\backslash
\{e\}\right\}  \text{.} \label{LipHn}%
\end{equation}
Of course, the sequence $\left(  H_{n}^{\prime}\right)  _{n\in\mathbb{N}}$ of
closed subspaces of $G$ converge to $G$ for $\mathfrak{H}(l)$, so we set
$H_{\infty}^{\prime}=G$. We observe, with obvious notations, that
$L_{H_{\infty}^{\prime}}=L_{\infty}$. Then the simple observation $L_{n}\leq
L_{H_{n}^{\prime}}\leq L_{n}^{\prime}$ lead to the following:

\begin{cor}
\label{LHM}For all $a\in V_{\varepsilon}$ we have $\lim_{n\rightarrow\infty
}L_{H_{n}^{\prime}}(a)=L_{\infty}(a)$.
\end{cor}

We can now conclude this section, by using \cite[Lemma 9.3]{Rieffel00}:

\begin{cor}
\label{Vcont}We have:%
\[
\lim_{n\rightarrow\infty}\operatorname*{dist}\nolimits_{q}\left(  \left(
\theta_{n}(V_{\varepsilon}),L_{n}\right)  ,\left(  V_{\varepsilon},L_{\infty
}\right)  \right)  =0\text{.}%
\]

\end{cor}

\subsection{Main Theorem}

We can now conclude this section with our main theorem, which is now an
immediate consequence of Proposition (\ref{finitedimapprox}) and Corollary
(\ref{Vcont}):

\begin{thm}
\label{MAIN}Let $G$ be a compact Abelian group endowed with a continuous
length function $l$. Let $\mathfrak{H}$ be the Hausdorff distance defined by
$l$ on the closed subsets of $G$. Let $(H_{n})_{n\in\mathbb{N}}$ be a sequence
of closed subgroups of $G$ converging to $G$ for $\mathfrak{H}$. Let
$\sigma_{\infty}$ be a skew bicharacter on the discrete Pontryagin dual
$\widehat{G}$ of $G$, and let for each $n\in\mathbb{N}$ a skew bicharacter
$\sigma_{n}$ of $\widehat{H_{n}}$ be given. We denote for all $n\in\mathbb{N}$
by $\sigma_{n}^{\prime}$ the lift of the bicharacter $\sigma_{n}$ to
$\widehat{G}$. We assume that the sequence $(\sigma_{n}^{\prime}%
)_{n\in\mathbb{N}}$ converges to $\sigma_{\infty}$ pointwise. Then:%
\[
\lim_{n\rightarrow\infty}\operatorname*{dist}\nolimits_{q}\left(  \left(
C^{\ast}\left(  \widehat{H_{n}},\sigma_{n}\right)  ,L_{n}\right)  ,\left(
C^{\ast}\left(  \widehat{G},\sigma_{\infty}\right)  ,L_{\infty}\right)
\right)  =0\text{.}%
\]

\end{thm}

\begin{pf}
Let $\varepsilon>0$. By Corollary (\ref{Vcont}), there exists $N\in\mathbb{N}$
such that $\operatorname*{dist}\nolimits_{q}\left(  \left(  \theta
_{n}(V_{\varepsilon}),L_{n}\right)  ,\left(  V_{\varepsilon},L_{\infty
}\right)  \right)  \leq\frac{1}{3}\varepsilon$ for all $n\geq N$. By the
triangle inequality property of $\operatorname*{dist}_{q}$ and by Proposition
(\ref{finitedimapprox}), we conclude that for all $n\geq\max\{N,N_{\varepsilon
}\}$ we have
\[
\operatorname*{dist}\nolimits_{q}\left(  \left(  C^{\ast}\left(
\widehat{H_{n}},\sigma_{n}\right)  ,L_{n}\right)  ,\left(  C^{\ast}\left(
\widehat{G},\sigma_{\infty}\right)  ,L_{\infty}\right)  \right)
\leq\varepsilon\text{.}%
\]This concludes our proof. \qed \end{pf}

Note that, using Corollary (\ref{LHM}), Theorem (\ref{MAIN}) remains true if
we replace the Lip-norms $L_{n}$ by the Lip-norms $L_{n}^{\prime}$.

\bigskip We shall apply Theorem (\ref{MAIN}) to obtain finite dimensional
approximations of the quantum tori.

\subsection{Finite Dimensional Approximations for the quantum Tori}

Our fundamental example (\ref{fundamentalexample}) satisfies Theorem
(\ref{MAIN}) and thus the conclusion announced in the introduction holds. We
wish to give a more specific application of Theorem (\ref{MAIN})\ where we
approximate some quantum tori by full matrix algebras. We use the notations of
example (\ref{fundamentalexample}):

\begin{prop}
\label{even}We assume $d$ is even. Suppose that $S_{\infty}$ is a block
\ diagonal matrix whose diagonal blocks are of the form $\psi_{j}\Lambda$ for
$j=1,\ldots,d/2$ where $\psi_{j}\in\lbrack0,1[$ and $\Lambda$ is the matrix%
\[
\left[
\begin{array}
[c]{cc}%
0 & -1\\
1 & 0
\end{array}
\right]  \text{.}%
\]
Let $\left(  p_{n}\right)  _{n\in\mathbb{N}}$ be the increasing sequence of
prime numbers in $\mathbb{N}\backslash\{0,1,2\}$. There exists a sequence
$(S_{n})_{n\in\mathbb{N}}$ of $d\times d$ antisymmetric matrices converging to
$S_{\infty}$ in operator norm, such that for all $n\in\mathbb{N}$ the
bicharacter $\sigma\lbrack S_{n}]$ of $\mathbb{Z}^{d}$ is also a bicharacter
of $\mathbb{Z}_{(p_{n},\ldots,p_{n})}^{d}$ and $C^{\ast}\left(  \mathbb{Z}%
_{(p_{n},\ldots,p_{n})}^{d},\sigma\lbrack S_{n}]\right)  $ is a full matrix
algebra. By Theorem (\ref{MAIN}), it follows that:%
\[
\left(  C^{\ast}\left(  \mathbb{Z}_{k_{n}}^{d},\sigma\lbrack S_{n}]\right)
,L_{n}\right)  \underset{n\rightarrow\infty}{\overset{\operatorname*{dist}%
\nolimits_{q}}{\longrightarrow}}\left(  C^{\ast}\left(  \mathbb{Z}^{d}%
,\sigma\lbrack S_{\infty}]\right)  ,L_{\infty}\right)  \text{.}%
\]

\end{prop}

\begin{pf}
We only sketch this easy argument (see \cite{Latremoliere04} for details). Let
$n\in\mathbb{N}$. Let $j\in\{1,\ldots,d/2\}$. We define $\psi_{j,n}$ by%
\[
\psi_{j,n}=\inf\left\{  m:m=1,\ldots,p_{n}-1\text{ and }\psi_{j}\leq\frac
{m}{p_{n}}\right\}  \text{.}%
\]
Since $\left\vert \psi_{j}-\frac{1}{p_{n}}\psi_{j,n}\right\vert \leq\frac
{1}{p_{n}}$ we conclude that $\lim_{n\rightarrow\infty}\frac{1}{p_{n}}%
\psi_{j,n}=\psi_{j}$. Let us define the matrix%
\[
S_{n}=\frac{1}{p_{n}}\left[
\begin{array}
[c]{cccc}%
\psi_{1,n}\Lambda & 0 & \cdots & 0\\
0 & \psi_{2,n}\Lambda & \ddots & \vdots\\
\vdots & \ddots & \ddots & 0\\
0 & \cdots & 0 & \psi_{d/2,n}\Lambda
\end{array}
\right]  \text{.}%
\]
We observe that $\sigma\lbrack S_{n}]$ defines a skew bicharacter on
$\mathbb{Z}_{(p_{n},\ldots,p_{n})}^{d}$. Now, the finite dimensional vector
space $C^{\ast}\left(  \mathbb{Z}_{(p_{n},\ldots,p_{n})}^{d},\sigma\lbrack
S_{n}]\right)  $ is spanned by the unitaries $\delta_{z}$ for $z\in
\mathbb{Z}_{(p_{n},\ldots,p_{n})}^{d}$ and $\delta_{z}$ the Dirac function at
$z$. Moreover $\delta_{z}\delta_{z^{\prime}}=\sigma\lbrack S_{n}](z,z^{\prime
})^{2}\delta_{z^{\prime}}\delta_{z}$, and the center of $C^{\ast}\left(
\mathbb{Z}_{(p_{n},\ldots,p_{n})}^{d},\sigma\lbrack S_{n}]\right)  $ is
spanned by those $\delta_{y}$ which commute with $\{\delta_{z}:z\in
\mathbb{Z}_{(p_{n},\ldots,p_{n})}^{d}\}$ --- in other words, the center of
$C^{\ast}\left(  \mathbb{Z}_{(p_{n},\ldots,p_{n})}^{d},\sigma\lbrack
S_{n}]\right)  $ is spanned by $\{\delta_{y}:\sigma\lbrack S_{n}%
](y,.)^{2}=1\}$. Yet $\sigma\lbrack S_{n}](y,.)^{2}=1$ if and only if
$2\psi_{\left\lceil j/2\right\rceil }y_{j}/p_{n}\in\mathbb{Z}$ for all
$j=1,\ldots,d$, where $y=(y_{1},\ldots,y_{j})$ and $\left\lceil
j/2\right\rceil $ is $\max\{j/2,(j+1)/2\}$. Since $0<\psi_{\left\lceil
j/2\right\rceil }<p_{n}$, and $p_{n}$ is a prime number strictly greater than
2, we conclude that $y_{j}$ must be divisible by $p_{n}$ in $\mathbb{Z}$, or
more precisely $y_{j}=0[p_{n}]$ for all $j=1,\ldots,d$. Thus with our choice
of $S_{n}$, the center of the finite dimensional C*-algebra $C^{\ast}\left(
\mathbb{Z}_{(p_{n},\ldots,p_{n})}^{d},\sigma\lbrack S_{n}]\right)  $ is
reduced to $\mathbb{C}1$ and therefore $C^{\ast}\left(  \mathbb{Z}%
_{(p_{n},\ldots,p_{n})}^{d},\sigma\lbrack S_{n}]\right)  $ is a full matrix algebra. \qed \end{pf}

As an immediate application, we get the following two corollaries:

\begin{cor}
Let $d$ be even and $d\geq2$. There exists a sequence $(A_{n},L_{n}%
)_{n\in\mathbb{N}}$ of compact quantum metric spaces such that for all
$n\in\mathbb{N}$ we have $A_{n}$ is a full matrix algebra, and
\[
\lim_{n\rightarrow0}\operatorname*{dist}\nolimits_{q}\left(  \left(
A_{n},L_{n}\right)  ,\left(  C\mathbb{(T}^{d}),L_{\infty}\right)  \right)
=0\text{.}%
\]

\end{cor}

\begin{pf}
Set $\psi_{0}=\ldots=\psi_{d/2-1}=0$ in Proposition (\ref{even}). \qed \end{pf}

\begin{cor}
\label{even1}Let $d=2$. Then for any $\sigma\in A(\mathbb{Z}^{2})$, there
exists a sequence $(A_{n},L_{n})_{n\in\mathbb{N}}$ of compact quantum metric
spaces such that for all $n\in\mathbb{N}$ we have $A_{n}$ is a full matrix
algebra, and
\[
\lim_{n\rightarrow0}\operatorname*{dist}\nolimits_{q}\left(  \left(
A_{n},L_{n}\right)  ,\left(  C\mathbb{(Z}^{2},\sigma),L_{\infty}\right)
\right)  =0\text{.}%
\]

\end{cor}

\begin{pf}
All skew bicharacters of $\mathbb{Z}^{2}$ are of the form $\sigma\lbrack
\psi\Lambda]$ for some $\psi\in\lbrack0,1[.$ We conclude by applying
Proposition (\ref{even}). \qed \end{pf}

We invite the interested reader to compare our construction in Corollary
(\ref{even1})\ and Proposition (\ref{even}) to the construction in
\cite{PimVoi80a}, where a sequence of direct sums of two full matrix algebras
has an inductive limit which is an AF algebra in which $C^{\ast}%
(\mathbb{Z}^{2},\sigma)$ embeds. This construction can also be found in
\cite[Secs. VI.3-VI.5, pp. 172-180]{Davidson}. The reader may check that each
summand in these finite dimensional algebras converge to $C^{\ast}%
(\mathbb{Z}^{2},\sigma)$ for $\operatorname*{dist}_{q}$ by Theorem
(\ref{MAIN}). We will return to the problem of approximating the tori of odd
dimension by full matrix algebras in the next section.

\section{Varying Length Functions and Dimensional Collapse}

We adopt a different framework in this section. We let $G$ be a compact group
with unit $e$, and $A$ be a unital C*-algebra whose norm is denoted by
$\left\Vert .\right\Vert _{A}$. We assume given a strongly continuous ergodic
action $\alpha$ of $G$ on $A$. The group $G$ is not necessarily Abelian, nor a
Lie group.

We let $(l_{n})_{n\in\mathbb{N}}$ be a sequence of continuous length functions
on $G$. \ We wish to deal with the more general case when the limit of a
sequence of length functions is not a length function on $G$, but factors into
a continuous length function on a quotient of $G$. Let $H$ be a normal
subgroup of $G$ and $K=G\backslash H$. By \cite[Theorem 5.26, p. 40]%
{Hewitt79}, the space $K$ is a topological group, and the canonical surjection
$g\in G\mapsto\lbrack g]\in K$ is continuous, so $K$ is compact.

Now, we suppose that the sequence $(l_{n})_{n\in\mathbb{N}}$ converges
uniformly on $G$ to a function $\widetilde{l_{\infty}}$ such that:%
\begin{equation}
\left\{
\begin{array}
[c]{c}%
\widetilde{l_{\infty}}(g)=0\text{ for }g\in H\text{,}\\
\widetilde{l_{\infty}}(g)>0\text{ for }g\not \in H\text{.}%
\end{array}
\right.  \label{Collapse0}%
\end{equation}

\begin{lem}
\label{ContL}If the sequence $(l_{n})_{n\in\mathbb{N}}$ satisfies the
conditions (\ref{Collapse0}), then there exists a continuous length function
$l_{\infty}$ on $K$ such that $l_{n}(g)\underset{n\rightarrow\infty
}{\longrightarrow}l_{\infty}\left(  [g]\right)  $ for all $g\in G$.
\end{lem}

\begin{pf}
Since $l_{n}$ is a length function on $G$ for all $n\in\mathbb{N}$, we have
that $\widetilde{l_{\infty}}(gg^{\prime})\leq\widetilde{l_{\infty}%
}(g)+\widetilde{l_{\infty}}(g^{\prime})$ and $\widetilde{l_{\infty}}%
(g^{-1})=\widetilde{l_{\infty}}(g)$ for all $g,g^{\prime}\in G$ by taking the
pointwise limit. Let $g\in G$ and $g^{\prime}\in H$. We observe that%
\begin{align*}
\widetilde{l_{\infty}}(g)  &  \leq\widetilde{l_{\infty}}(gg^{\prime
})+\widetilde{l_{\infty}}\left(  \left(  g^{\prime}\right)  ^{-1}\right)
=\widetilde{l_{\infty}}(gg^{\prime})\\
&  \leq\widetilde{l_{\infty}}(g)+\widetilde{l_{\infty}}(g^{\prime}%
)=\widetilde{l_{\infty}}(g)
\end{align*}
since $\widetilde{l_{\infty}}(g^{\prime})=\widetilde{l_{\infty}}\left(
\left(  g^{\prime}\right)  ^{-1}\right)  =0$ because $g^{\prime}\in H$ and
(\ref{Collapse0}) holds by assumption. We thus have shown that $\widetilde
{l_{\infty}}(g)=\widetilde{l_{\infty}}(gg^{\prime})$ for all $g^{\prime}\in
H$. We can thus define a function $l_{\infty}$ on $K$ by setting $l_{\infty
}([g])=\widetilde{l_{\infty}}(g)$, and by construction $\lim_{n\rightarrow
\infty}l_{n}(g)=l_{\infty}\left(  [g]\right)  $. We have, for all
$g,g^{\prime}\in G$:%
\[
l_{\infty}([g][g]^{\prime})=\widetilde{l_{\infty}}(gg^{\prime})\leq
\widetilde{l_{\infty}}(g)+\widetilde{l_{\infty}}(g^{\prime})=l_{\infty
}([g])+l_{\infty}([g^{\prime}])\text{,}%
\]
and similarly $l_{\infty}\left(  [g]^{-1}\right)  =l_{\infty}([g])$. Moreover,
$l_{\infty}([g])=0$ if and only if $\widetilde{l_{\infty}}(g)=0,$ which is
equivalent by (\ref{Collapse0}) to $g\in H$. Hence $l_{\infty}\left(  \left[
g\right]  \right)  =0$ if and only if $[g]$ is the identity $[e] $ of $K$.
Therefore $l_{\infty}$ is a length function on $K$. It remains to prove that
$l_{\infty}$ is continuous. Now, the canonical surjection $g\in G\mapsto
\lbrack g]\in K$ is continuous and open by \cite[Theorem 5.17, p.
37]{Hewitt79}. Therefore, $l_{\infty}$ is continuous on $K=[G]$ if and only if
$g\mapsto l_{\infty}\circ\lbrack g]$ is continuous on $G$. Yet by construction
$g\mapsto l_{\infty}\circ\lbrack g]$ is $\widetilde{l_{\infty}}$, which is
continuous on $G$, since it is the uniform limit of the sequence
$(l_{n})_{n\in\mathbb{N}}$ of continuous (length) functions. This concludes
our lemma. \qed \end{pf}

\begin{rem}
\bigskip We observe that, without assuming (\ref{Collapse0}), the zero set of
the pointwise limit $\widetilde{l_{\infty}}$ of the sequence $(l_{n}%
)_{n\in\mathbb{N}}$ of continuous length functions on $G$ is always a group.
It is also straightforward that this group is normal if there exists some
$N\in\mathbb{N}$ such that for all $n\geq N$, the length functions $l_{n}$ are
invariant by the inner automorphisms of $G$. If moreover the sequence
$(l_{n})_{n\in\mathbb{N}}$ converges uniformly, we can set $H$ to be this
normal subgroup of $G$ and by construction the hypotheses of Lemma
(\ref{ContL}) are then satisfied.
\end{rem}

For each $n\in\mathbb{N}$, we define on $K$ the quotient length functions
$l_{n}^{K}$ by setting for all $g\in G$:%
\[
l_{n}^{K}([g])=\inf\left\{  l_{n}(g^{\prime}):g^{\prime}\in G\text{,
}[g]=[g^{\prime}]\right\}  \text{.}%
\]

The following lemma will be useful:

\begin{lem}
\label{Lnk}For each $n\in\mathbb{N}$, the length function $l_{n}^{K}$ is
continuous on $K$. Moreover, the sequence $(l_{n}^{K})_{n\in\mathbb{N}}$
converges uniformly on $K$ to $l_{\infty}$.
\end{lem}

\begin{pf}
Let $n\in\mathbb{N}$. By definition, we have for all $[g]\in K$:%
\[
l_{n}^{K}([g])=\inf\left\{  l_{n}(gh):h\in H\right\}  \text{.}%
\]

Now, since $l_{n}$ is a length function, we have for all $h\in H$ and
$g^{\prime}\in G$:%
\[
\left\vert l_{n}(gh)-l_{n}(g^{\prime}h)\right\vert \leq l_{n}(g^{-1}g^{\prime
})
\]
and so the family of functions $\left(  g\mapsto l_{n}(gh)\right)  _{h\in H}$
is equicontinuous on the compact $G$. Their infimum is therefore continuous as
well, so $g\in G\mapsto l_{n}^{K}([g])$ is continuous. Since $g\in
G\mapsto\lbrack g]\in K$ is a continuous open surjection, we deduce that
$l_{n}^{K}$ is continuous on $K$.

Now let us fix $g\in G$. Since $H$ is compact and $h\mapsto l_{n}(gh)$ is
continuous on $H$, there exists $h_{n}\in H$ such that $\inf\left\{
l_{n}(gh):h\in H\right\}  =l_{n}(gh_{n})$. We then have $\left\vert l_{n}%
^{K}([g])-l_{\infty}([g])\right\vert =\left\vert l_{n}(gh_{n})-\widetilde
{l_{\infty}}(gh_{n})\right\vert $ since we have shown in Lemma (\ref{ContL})
that $\widetilde{l_{\infty}}$ is constant on the coset $gH$ in $G$. Therefore:%
\begin{equation}
\sup_{\lbrack g]\in K}\left\vert l_{n}^{K}([g])-l_{\infty}([g])\right\vert
\leq\sup_{g\in G}\left\vert l_{n}(g)-\widetilde{l_{\infty}}(g)\right\vert
\text{.} \label{Lnk0}%
\end{equation}
By definition, $\widetilde{l_{\infty}}$ is the uniform limit on $G$ of the
sequence $(l_{n})_{n\in\mathbb{N}}$ so $(l_{n}^{K})_{n\in\mathbb{N}}$
converges uniformly on $K$ to $l_{\infty}$ by (\ref{Lnk0}). \qed \end{pf}

We recall that, given an ergodic strongly continuous action $\alpha$ of any
compact group $\Sigma$ of unit $e$ with a length function $\rho$ on a unital
C*-algebra $A$ of norm $\left\Vert .\right\Vert _{A}$, Rieffel proved in
\cite{Rieffel98a} that the seminorm defined for all $a\in A$ by%
\begin{equation}
L(a)=\sup\left\{  \rho(g)^{-1}\left\Vert a-\alpha_{g}(a)\right\Vert
:g\in\Sigma\backslash\{e\}\right\}  \label{lipnormbyaction}%
\end{equation}
is a Lip-norm on $A$.

We are now ready to prove the main result of this section:

\begin{thm}
\label{Collapse}Let $A_{K}$ be the fixed-point C*-algebra of the action
$\alpha$ restricted to $H$. Let $(l_{n})_{n\in\mathbb{N}}$ be a sequence of
continuous length functions on $G$ satisfying the conditions (\ref{Collapse0}%
). Moreover, we assume that%
\begin{equation}
\lim_{n\rightarrow\infty}\sup_{[g]\in K\backslash\{[e]\}}\left\vert
\frac{l_{\infty}([g])}{l_{n}^{K}([g])}-1\right\vert =0 \label{Collapse00}%
\end{equation}
and that for all $n\in\mathbb{N}$ we have $l_{n}(g^{-1}g^{\prime
}g)=l(g^{\prime})$ for all $g,g^{\prime}\in G$.

For all $n\in\mathbb{N}$ we define on $A$ the Lip-norm $L[n]$ from $\alpha$
and $l_{n}$ as in (\ref{lipnormbyaction}). Now, let $\alpha^{K}$ be the
quotient action of $K$ on $A_{K}$ induced by $\alpha$. The action $\alpha^{K}$
is strongly continuous and ergodic. We can therefore define the Lip norm
$L[\infty]$ on $A_{K}$ obtained from the action $\alpha^{K}$ of $K$ on $A_{K}$
and the continuous length function $l_{\infty}$ of Lemma\ (\ref{ContL}). We
have:%
\[
\lim_{n\rightarrow\infty}\operatorname*{dist}\nolimits_{q}\left(  \left(
A,L[n]\right)  ,\left(  A_{K},L[\infty]\right)  \right)  =0\text{.}%
\]

\end{thm}

\begin{pf}
We start by remarking that $\alpha$ is a strongly continuous ergodic action on
$A_{K}$ which coincides with an action $\alpha^{K}$ of $K$ on $A_{K}$. Indeed,
let $a\in A_{K}$, $g\in G$ and $h\in H$. Then $\alpha_{gh}(a)=\alpha
_{g}(\alpha_{h}(a))$ but by definition of $A_{K}$ we have $\alpha_{h}(a)=a$.
Hence we can define $\alpha_{\lbrack g]}^{K}(a)$ by $\alpha_{g}(a)$ for any
$g\in G$ and for all $a\in A_{K}$, where as usual $[g]$ is the class of $g$ in
$K$. Moreover we have:%
\begin{align*}
\alpha_{h}(\alpha_{\lbrack g]}^{K}(a))  &  =\alpha_{hg}(a)=\alpha_{g}\left(
\alpha_{g^{-1}hg}(a)\right) \\
&  =\alpha_{g}(a)=\alpha_{\lbrack g]}^{K}(a)\text{,}%
\end{align*}
since $h\in H$ and so $g^{-1}hg\in H$ since $H$ is a normal subgroup of $G$.
Therefore, since $h\in H$ was chosen arbitrarily, $\alpha_{\lbrack g]}%
^{K}(a)\in A_{K}$ by definition. Hence, $\alpha^{K}$ is an action of $K$ on
$A_{K}$. Since $\alpha$ is ergodic, so is $\alpha^{K}$. Since $g\in
G\mapsto\lbrack g]\in K$ is a continuous open surjection of $G$ onto $K$, the
strong continuity of $\alpha$ on $G$ implies the strong continuity of
$\alpha^{K}$ on $K$.

Thus, the action $\alpha^{K}$ and the continuous length function $l_{\infty}$
indeed define a Lip-norm $L[\infty]$ on $A_{K}$ by setting, for all $a\in
A_{K}$:%
\[
L[\infty](a)=\sup\left\{  \frac{\left\Vert a-\alpha_{\lbrack g]}%
^{K}(a)\right\Vert _{A}}{l_{\infty}([g])}:[g]\in K\backslash\{[e]\}\right\}
\text{.}%
\]

Now we also observe that for all $n\in\mathbb{N}$, and most importantly for
all $a\in A_{K}$:%
\begin{align}
L[n](a)  &  =\sup\left\{  \frac{\left\Vert a-\alpha_{g}(a)\right\Vert _{A}%
}{l_{n}(g)}:g\in G\backslash\{e\}\right\} \nonumber\\
&  =\sup\left\{  \frac{\left\Vert a-\alpha_{\lbrack g]}^{K}(a)\right\Vert
_{A}}{l_{n}^{K}([g])}:[g]\in K\backslash\{[e]\}\right\}  \text{.} \label{lnk}%
\end{align}

In other words, for each $n\in\mathbb{N}$, the restriction of $L[n]$ to
$A_{K}$ agrees with the Lip-norm $L^{K}[n]$ obtained from the action
$\alpha^{K}$ of $K$ and the continuous length function $l_{n}^{K}$ of $K$. We
shall only use the notation $L[n]$. Also, we observe that for any
$n\in\mathbb{N}$, the Lip-norm $L[n]$ is the supremum of a family of
continuous functions on $A$, so it is lower-semicontinuous on $A$. Similarly,
$L[\infty]$ is lower-semicontinuous on $A_{K}$.

We now turn to the main argument of our theorem. We define the conditional
expectation $E$ on $A$ by setting for all $a\in A$:
\[
E(a)=\int_{h\in H}\alpha_{h}(a)d\lambda_{H}(h)
\]
where $\lambda_{H}$ is the Haar probability measure on $H$. A standard
argument shows that $E$ is indeed a conditional expectation whose range is the
fixed point C*--algebra $A_{K}$ of the action $\alpha$ restricted to\ $H$.
Now, we can compute:%
\begin{align}
\left\Vert a-E(a)\right\Vert _{A}  &  =\left\Vert \int_{H}\left(  a-\alpha
_{h}(a)\right)  d\lambda_{H}(h)\right\Vert _{A}\leq\int_{H}\left\Vert
a-\alpha_{h}(a)\right\Vert _{A}d\lambda_{H}(h)\nonumber\\
&  \leq\int_{H}L[n](a)l_{n}(h)d\lambda_{H}(h)=L[n](a)\int_{H}l_{n}%
(h)d\lambda_{H}(h) \label{collapse0.1}%
\end{align}
for all $n\in\mathbb{N}$. We used in the first equality the fact that
$\lambda_{H}$ is a probability measure.

Let $\varepsilon>0$ be given. Now, by assumption, $(l_{n})_{n\in\mathbb{N}}$
converges uniformly to 0 on $H$, hence $\lim_{n\rightarrow\infty}\int_{H}%
l_{n}(h)d\lambda_{H}(h)=0$. Therefore there exists $N_{0}\in\mathbb{N}$ such
that for all $n\geq N_{0}$ we have $\int_{H}l_{n}(h)d\lambda_{H}(h)\leq
\frac{1}{4}\varepsilon$. Hence by (\ref{collapse0.1}) we have:%
\begin{equation}
\left\Vert a-E(a)\right\Vert _{A}\leq\frac{1}{4}\varepsilon L[n](a)\text{.}
\label{collapse0.2}%
\end{equation}

Since $L[n]$ is a lower semicontinuous seminorm, we have that
\[
L[n]\left(  \int_{H}\alpha_{h}(a)d\lambda_{H}(h)\right)  \leq\int
_{H}L[n](\alpha_{h}(a))d\lambda_{H}(h)\text{.}%
\]
Now we have that $L[n](\alpha_{h}(a))=L[n](a)$ since the inner automorphisms
of $G$ preserve $l_{n}$. Hence we have:
\[
L[n]\left(  \int_{H}\alpha_{h}(a)d\lambda_{H}(h)\right)  \leq\left(  \int
_{H}d\lambda_{H}\right)  L[n](a)=L[n](a)
\]
as $\lambda_{H}$ is a probability measure on $H$. In other words, for all
$a\in A$ we have:%
\begin{equation}
L[n](E(a))\leq L[n](a)\text{.} \label{collapse0.3}%
\end{equation}

Hence, using \cite[Proposition 8.5]{Rieffel00} with (\ref{collapse0.2}) and
(\ref{collapse0.3}), we have, for all $n\geq N_{0}$:
\begin{equation}
\operatorname*{dist}\nolimits_{q}\left(  \left(  A,L[n]\right)  ,\left(
A_{K},L[n]\right)  \right)  \leq\frac{1}{4}\varepsilon. \label{collapse1}%
\end{equation}

We now use a similar method to that of the proof of Proposition
(\ref{finitedimapprox}). Let us endow $K$ with its unique probability Haar
measure $\lambda_{K}$. For any $p\in\mathbb{N}\backslash\{0\}$, we denote by
$\varphi_{p}$ the function given by Lemma (\ref{fejer}) such that $\int
_{G}l_{\infty}\varphi_{p}d\lambda<p^{-1}$. Let $p\in\mathbb{N}\backslash\{0\}$
such that
\begin{equation}
\int_{K}l_{\infty}(g)\varphi_{p}(g)d\lambda_{K}(g)\leq\frac{1}{8}%
\varepsilon\text{.} \label{collapse2}%
\end{equation}
We define the operator $\alpha_{\varphi_{p}}^{K}$ on $A_{K}$ by setting, for
all $a\in A_{K}$:%
\[
\alpha_{\varphi_{p}}^{K}(a)=\int_{K}\varphi_{p}(g)\alpha_{g^{-1}}%
^{K}(a)d\lambda_{K}(g)\text{.}%
\]
We denote by $V_{p}$ the image by $\alpha_{\varphi_{p}}^{K}$ of the
self-adjoint part of $A_{K}$. It is an easy corollary of \cite[sec. 8, pp.
37--39]{Rieffel00} that $V_{p}$ is an finite dimensional order-unit subspace
of $A_{K}$ such that $L[\infty](a)<\infty$ for all $a\in V_{p}$ (specifically,
$A_{K}$ is the direct sum of the isotopic components of $A$ for the
irreducible representations whose characters appear in $\varphi_{p}$). Now, as
in the proof of Proposition (\ref{finitedimapprox}), we have, for all
$n\in\overline{\mathbb{N}}$:%
\begin{equation}
\left\Vert a-\alpha_{\varphi_{p}}^{K}(a)\right\Vert _{A}\leq\left(  \int
_{K}l_{n}^{K}(g)\varphi_{p}(g)d\lambda_{K}(g)\right)  L[n](a)\text{,}
\label{movelength0.1}%
\end{equation}
where we denoted $l_{\infty}$ by $l_{\infty}^{K}$ to ease notations. Now, by
Lemma (\ref{Lnk}), the sequence $\left(  l_{n}^{K}\right)  _{n\in\mathbb{N}}$
converges uniformly to $l_{\infty}$ on $K$. Therefore, there exists $N_{1}%
\in\mathbb{N}$ such that, for all $n\geq N_{1}$ we have:%
\begin{equation}
\int_{K}l_{n}^{K}(g)\varphi_{p}(g)d\lambda_{K}(g)\leq\frac{1}{4}%
\varepsilon\label{movelength0.2}%
\end{equation}
using (\ref{collapse2}). We conclude from (\ref{movelength0.1}) and
(\ref{movelength0.2}) that for all $a\in A_{K}$ and all $n\in\overline
{\mathbb{N}}$ such that $n\geq N_{1}$ we have:%
\begin{equation}
\left\Vert a-\alpha_{\varphi_{p}}^{K}(a)\right\Vert _{A}\leq\frac{1}%
{4}\varepsilon L[n](a)\text{.} \label{movelength0.3}%
\end{equation}

Of course, with the same argument as used in the proof of Proposition
(\ref{finitedimapprox}) and also used to deduce (\ref{collapse0.3}), we have,
for all $n\in\overline{\mathbb{N}}$ such that $n\geq N_{1}$ and for all $a\in
A_{K}$:%
\begin{equation}
L[n](\alpha_{\varphi_{p}}^{K}(a))\leq L[n](a)\text{.} \label{movelength0.4}%
\end{equation}

By construction, the range of $\alpha_{\varphi_{p}}^{K}$ restricted to the
self-adjoint part of $A_{K}$ is the order-unit subspace $V_{p}$ of $A_{K}$.
Therefore, by \cite[Proposition 8.5]{Rieffel00} applied to
(\ref{movelength0.3}) and (\ref{movelength0.4}), we have that for all
$n\in\overline{\mathbb{N}}$ such that $n\geq N_{1}$:%
\begin{equation}
\operatorname*{dist}\nolimits_{q}\left(  \left(  \left(  A_{K},L[n]\right)
\right)  ,\left(  V_{p},L[n]\right)  \right)  \leq\frac{\varepsilon}%
{4}\text{.} \label{MoveLength1}%
\end{equation}
Our situation is simpler than the setting of Proposition
(\ref{finitedimapprox}) in which $V_{p}$ was viewed as a subspace of many
different C*--algebras.

Now, let $a\in V_{p}\subseteq A_{K}.$ Let $g\in K\backslash\{[e]\}$. We have:%
\begin{align*}
\left\Vert \alpha_{g}^{K}(a)-a\right\Vert _{A}\left\vert \frac{1}{l_{n}%
^{K}(g)}-\frac{1}{l_{\infty}(g)}\right\vert  &  \leq L[\infty](a)l_{\infty
}(g)\left\vert \frac{1}{l_{n}^{K}(g)}-\frac{1}{l_{\infty}(g)}\right\vert \\
&  =L[\infty](a)\left\vert \frac{l_{\infty}(g)}{l_{n}^{K}(g)}-1\right\vert
\text{.}%
\end{align*}
Let $\varepsilon^{\prime}>0$ be given. By assumption we have $L[\infty
](a)<\infty$ and%
\[
\lim_{n\rightarrow\infty}\sup_{g\in K\backslash\{[e]\}}\left\vert
\frac{l_{\infty}(g)}{l_{n}^{K}(g)}-1\right\vert =0.
\]
Therefore there exists $N_{2}\in\mathbb{N}$ such that for all $n\geq N_{2}$
and for all $g\in K\backslash\{[e]\}$ we have $\left\vert \frac{l_{\infty}%
(g)}{l_{n}^{K}(g)}-1\right\vert \leq\varepsilon^{\prime}$ and thus:
\[
\left\Vert \alpha_{g}^{K}(a)-a\right\Vert _{A}\left\vert \frac{1}{l_{n}%
^{K}(g)}-\frac{1}{l_{\infty}(g)}\right\vert \leq\varepsilon^{\prime}%
L[\infty](a)\text{,}%
\]
which implies that $\left\vert L[n](a)-L[\infty](a)\right\vert \leq
\varepsilon^{\prime}\left(  L[\infty](a)\right)  $, using our observation
(\ref{lnk}). Hence $\lim_{n\rightarrow\infty}L[n](a)=L[\infty](a)$.

Summarizing, on the finite dimensional order-unit space $V_{p}$ we have a
continuous field of Lip-norms $\left(  L[n]\right)  _{n\in\overline
{\mathbb{N}}}$, and so by \cite[theorem 11.2]{Rieffel00} we conclude there
exists $N_{3}\in\mathbb{N}$ such that for all $n\geq N_{3}$ we have:%
\begin{equation}
\operatorname*{dist}\nolimits_{q}\left(  \left(  V_{p},L[n]\right)
,(V_{p},L[\infty])\right)  \leq\frac{1}{4}\varepsilon. \label{MoveLength2}%
\end{equation}
Hence by the triangle inequality, using (\ref{MoveLength1}) and
(\ref{MoveLength2}), we have for all $n\geq\max\{N_{1},N_{3}\}$:
\begin{equation}
\operatorname*{dist}\nolimits_{q}\left(  \left(  \left(  A_{K},L[n]\right)
\right)  ,\left(  \left(  A_{K},L[\infty]\right)  \right)  \right)  \leq
\frac{3}{4}\varepsilon\text{.} \label{MoveLength3}%
\end{equation}

By the triangle inequality, using (\ref{MoveLength3}) and (\ref{collapse1}),
we conclude that for all $n\geq N_{4}$ where $N_{4}=\max\{N_{0},N_{1,}N_{3}%
\}$:%
\[
\operatorname*{dist}\nolimits_{q}\left(  \left(  \left(  A,L[n]\right)
\right)  ,\left(  \left(  A_{K},L[\infty]\right)  \right)  \right)
\leq\varepsilon\text{.}%
\]

This concludes the proof of this Theorem. \qed \end{pf}

\begin{rem}
By Lemma (\ref{Lnk}), for each $n\in\mathbb{N}$ the length function $l_{n}%
^{K}$ is continuous on $K$, so it can be used to define a Lip-norm $L^{K}[n]$
on $A_{K}$ via the strongly continuous ergodic action $\alpha^{K}$. The proof
of Theorem (\ref{Collapse}) shows in particular that the field of Lip-norms
$n\in\overline{\mathbb{N}}\mapsto L^{K}[n]$ (where $L^{K}[\infty]=L[\infty]$)
is continuous on $A_{K}$.
\end{rem}

\bigskip We obtain in particular the following corollary, when $H=\{e\}$:

\begin{cor}
Let $(l_{n})_{n\in\mathbb{N}}$ be a sequence of continuous length functions on
$G$ uniformly converging on $G$ to a function $l_{\infty}$, such that%
\[
\lim_{n\rightarrow\infty}\sup_{g\in G\backslash\{e\}}\left\vert \frac
{l_{\infty}(g)}{l_{n}(g)}-1\right\vert =0
\]
and $l_{\infty}(g)=0$ only if $g=e$. Then $l_{\infty}$ is a continuous length
function on $G$. If we denote by $L[n]$ the Lip-norm defined on $A$ by
$\alpha$ and $l_{n}$ for all $n\in\overline{\mathbb{N}}$ then
\[
\lim_{n\rightarrow\infty}\operatorname*{dist}\nolimits_{q}\left(  \left(
A,L[n]\right)  ,\left(  A,L[\infty]\right)  \right)  =0
\]

\end{cor}

\begin{pf}
Setting $H=\{e\}$, we can first use Lemma (\ref{ContL}) to conclude that
$l_{\infty}$ is indeed a continuous length function on $K=G$. Then we simply
apply Theorem (\ref{Collapse}) to complete the proof of our corollary. \qed \end{pf}

\bigskip Theorem (\ref{Collapse}) is not vacuous, as there exists sequences of
continuous length functions satisfying its hypotheses. We now focus on our
fundamental example (\ref{fundamentalexample}). The following easy lemma
proposes a construction of such sequences when $G=\mathbb{U}_{k}^{d}$ and
$H=\{1\}\times\mathbb{U}_{k^{\prime\prime}}^{d-\delta}$ for any $k\in
\mathbb{N}^{d}$, any $\delta\in\{0,\ldots,d\}$, and where $k^{\prime\prime}%
\in\mathbb{N}^{d-\delta}$ and $k^{\prime}\in\mathbb{N}^{\delta}$ are uniquely
defined by $k=(k^{\prime},k^{\prime\prime})$. With these notations,
$K=\mathbb{U}_{k^{\prime}}^{\delta}$. Observe that we can also write
$\mathbb{U}_{k}^{d}=\mathbb{U}_{k^{\prime}}^{\delta}\times\mathbb{U}%
_{k^{\prime\prime}}^{d-\delta}$. We have:

\begin{lem}
\label{Collapse3}Let $l$ be a continuous length function on $\mathbb{U}%
_{k}^{d}$. We set, for $n\in\mathbb{N}$:%
\[
l_{n}(\omega_{1},\omega_{2})=\frac{1}{n+1}l(\omega_{1},\omega_{2})+\left(
1-\frac{1}{n+1}\right)  l(\omega_{1},1)
\]
where $\omega_{1}\in\mathbb{U}_{k^{\prime}}^{\delta}$ and $\omega_{2}%
\in\mathbb{U}_{k^{\prime\prime}}^{d-\delta}$. Then the sequence $(l_{n}%
)_{n\in\mathbb{N}}$ satisfies the hypothesis of Theorem (\ref{Collapse}),
namely (\ref{Collapse0}) and (\ref{Collapse00}).Moreover the length function
$l_{\infty}$ on $\mathbb{U}_{k^{\prime}}^{\delta}$ is given by $\omega
\in\mathbb{U}_{k^{\prime}}^{\delta}\longmapsto l(\omega,1)$.
\end{lem}

\begin{pf}
This proof is straightforward. We refer to \cite{Latremoliere04}. \qed \end{pf}

An interesting remark is that we can use Theorem (\ref{Collapse}) together
with Lemma (\ref{Collapse3}) to prove that there exist sequences of quantum
tori converging for $\operatorname*{dist}_{q}$ to a circle ($\delta=1$) or
even a point ($\delta=0$).

\bigskip Now, we have the following obvious corollary from Theorem
(\ref{Collapse}) and Proposition (\ref{even}):

\begin{cor}
\label{odd}Let $d\in2\mathbb{N}+1$. There exists a sequence $(A_{n},L_{A_{n}%
})_{n\in\mathbb{N}}$ of compact quantum metric spaces, where $A_{n}$ is a full
matrix algebra for each $n\in\mathbb{N}$, such that $(A_{n},L_{A_{n}}%
)_{n\in\mathbb{N}}$ converges to the ordinary torus $\left(  C(\mathbb{T}%
^{d}),L_{\infty}\right)  =\left(  C^{\ast}(\mathbb{Z}^{d}),L_{\infty}\right)
$ for $\operatorname*{dist}_{q}$.
\end{cor}

\begin{pf}
This is a simple diagonal argument. Let $n\in\mathbb{N}$. Let $l_{m}$ be given
as in Lemma\ (\ref{Collapse3}) for any $m\in\mathbb{N}$. By Proposition
(\ref{even})\ there exists a prime number $q_{n,m}$ and a skew bicharacter
$\sigma_{n,m}$ of $\mathbb{Z}_{(q_{n,m},\ldots,q_{n,m})}^{d+1}$ such that
$A_{n,m}=C^{\ast}\left(  \mathbb{Z}_{(q_{n,m},\ldots,q_{n,m})}^{d+1}%
,\sigma_{n,m}\right)  $ is a full matrix algebra and
\[
\operatorname*{dist}\nolimits_{q}\left(  (A_{n,m},L_{A_{n,m}}),(C(\mathbb{T}%
^{d+1}),L[l_{m}])\right)  \leq(2n)^{-1}%
\]
where $L_{A_{n,m}}$ is the lip-norm defined on $C^{\ast}\left(  \mathbb{Z}%
_{(q_{n,m},\ldots,q_{n,m})}^{d+1},\sigma_{n,m}\right)  $ by the dual action of
$\mathbb{U}_{(q_{n,m},\ldots,q_{n,m})}^{d+1}$ and the length $l_{m}$. Now,
using Theorem (\ref{Collapse}), we can pick $m$ large enough so that
\[
\operatorname*{dist}\nolimits_{q}\left(  (C(\mathbb{T}^{d}%
),L[l]),(C(\mathbb{T}^{d+1}),L[l_{m}])\right)  \leq(2n)^{-1}\text{.}%
\]
Hence for this $m$ we have $\operatorname*{dist}\nolimits_{q}\left(
(A_{n,m},L_{A_{n,m}}),(C(\mathbb{T}),L[l])\right)  \leq n^{-1}$. \qed \end{pf}

\section{Annex: Any compact metric space can be approximated in
Gro\-mov--Hausdorff distance by matrix algebras}

One could be tempted to ask the question: when is any compact metric space the
limit of full matrix algebras for $\operatorname{dist}_{q}$? Surprisingly, the
answer is all of them. However, the construction we will propose is blind to
any symmetry of the limit space, and in particular, it does not prove our
theorem (\ref{MAIN}) even when the limit is $C^{\ast}(\mathbb{Z}%
^{d})=C(\mathbb{T}^{d})$, as the Lip-norms on the matrix algebras are not
coming from any action of a product of cyclic groups. Of course, the following
result does not tell us anything about approximations of the quantum tori in general.

Let $\left(  X=\{x_{1},\ldots,x_{n}\},d\right)  $ be a metric space of
cardinality $n\in\mathbb{N}$ whose distance $d$ induces a Lipschitz seminorm
$L$ on $C(X)$. Let $\varepsilon>0$. Also, define the following maps:
\begin{align*}
D_{n}  &  :f\in C(X)\longmapsto\left[
\begin{array}
[c]{cccc}%
f(x_{1}) & 0 & \cdots & 0\\
0 & f(x_{2}) & \ddots & \vdots\\
\vdots & \ddots & \ddots & 0\\
0 & \cdots & 0 & f(x_{n})
\end{array}
\right]  \text{,}\\
P_{n}  &  :A\in M_{n}\longmapsto\left(  x_{i}\mapsto A_{i,i}\right)  \text{,}%
\end{align*}
where for any matrix $A$, the complex number $A_{i,j}$ is the
$(i,j)^{\text{th}}$ entry of $A$. We set the seminorm:
\[
L_{o,n}(A)=\sup_{i\neq j}\left\vert A_{ij}\right\vert \text{.}%
\]

With these notations introduced, we can state our result:

\begin{prop}
\label{easycv}Let $\varepsilon>0$ be given. Set, for $n\geq1$:
\[
L_{\varepsilon}=\max\{L\circ P_{n},\frac{n^{2}-n}{2\varepsilon}L_{o,n}%
\}\text{.}%
\]
Then $L_{\varepsilon}$ is a Lip-norm on $M_{n}$, and moreover
\[
\operatorname{dist}_{q}((C(X),L),(M_{n},L_{\varepsilon}))\leq\varepsilon
\text{.}%
\]

\end{prop}

\begin{pf}
Proving that a seminorm such as $L_{\varepsilon}$ is a Lip-norm for a finite
dimensional algebra is a triviality, as all that is needed is that
$L_{\varepsilon}(A)=0$ if, and only if $A=\lambda1$ for some $\lambda
\in\mathbb{C}$. Now, if $L_{\varepsilon}(A)=0$ then $L_{o,n}(A)=0$ so $A$ is
diagonal, and since $L(P(A))=0$, $P(A)$ is a constant function, hence
$L_{\varepsilon}$ is a Lip-norm.

We now set, for $f\in C(X)$ and $A\in M_{n}$:
\begin{align*}
N(f,A)  &  =\frac{2}{\varepsilon}\left\Vert f-P_{n}(A)\right\Vert
_{C(X)}\text{,}\\
L_{N}(f,A)  &  =\max\{L(f),L_{\varepsilon}(A),N(f,A)\}\text{.}%
\end{align*}
$N$ is trivially a continuous seminorm over $C(X)\oplus M_{n}$, and we check:
\[
N(1,0)=N(0,1)=\frac{2}{\varepsilon}%
\]
and yet $N(1,1)=0$ and, obviously,
\[
N(f,D_{n}(f))=0=N(P_{n}(A),A)\text{.}%
\]

Now, this implies that $N$ is a bridge for $L$ and $L_{\varepsilon}$ as
defined in \cite{Rieffel01}, as indeed we check:
\[
\left\{
\begin{array}
[c]{c}%
L(P_{n}(A))\leq L_{\varepsilon}(A)\\
L_{\varepsilon}(D_{n}(f))=L(f)\text{.}%
\end{array}
\right.
\]
In other words, $L_{N}$ is a Lip-norm compatible with $(C(X),L)$ and
$(M_{n},L_{\varepsilon})$. It is therefore now sufficient to prove that the
Hausdorff distance between the state spaces of $C(X)$ and $M_{n}$ in the state
space of $C(X)\oplus M_{n}$ are $\varepsilon$-close. Let $(\mu,\nu)\in
S(C(X))\times S(M_{n})$, and let $(f,A)\in C(X)\oplus M_{n}$ be such that
$L_{N}(f,A)\leq1$. Then by definition of $N$, we have $\left\Vert
f-P_{n}(A)\right\Vert \leq\frac{1}{2}\varepsilon$, so, defining the linear
functional $\mu\circ P_{n}$ on $M_{n}$, we see that
\[
\left\vert \mu(f)-\mu(P_{n}(A))\right\vert \leq\frac{1}{2}\varepsilon\text{.}%
\]

Since $P_{n}$ is a positive map and $P_{n}$ maps the unit of $M_{n}$ to the
unit of $C(X)$, the linear map $\mu\circ P_{n}$ is a state of $M_{n}$. Thus
$S(C(X))$ is in an $\varepsilon$-neighborhood of $S(M_{n})$ in Hausdorff
distance for the dual metric of $L_{N}$.

On the other hand, $D$ being trivially a positive unit-preserving map, the map
$\nu\circ D$ is a state of $C(X)$, and:%
\begin{align*}
\left\vert \nu(D(f))-\nu(A)\right\vert  &  \leq\left\vert \nu(A-D\circ
P(A))\right\vert +\left\vert \nu(D(f)-D\circ P(A))\right\vert \\
&  \leq\left\Vert A-D\circ P(A)\right\Vert +\left\Vert f-P(A)\right\Vert \\
&  \leq(n^{2}-n)L_{o,n}(A)+\frac{1}{2}\varepsilon\leq\varepsilon\text{.}%
\end{align*} \qed \end{pf}

\begin{cor}
Any compact metric space $(X,d)$ is a limit point of the space of full matrix
algebras with Lip-norms, topologized by the quantum Gromov-Hausdorff distance.
\end{cor}

\begin{pf}
Let $\varepsilon>0$. Since $X$ is compact, there exists a finite subset
$X_{\varepsilon}$ of $X$ such that $X\subseteq\bigcup_{x\in X_{n}%
}B(x,\varepsilon)$ where $B(x,\varepsilon)$ is the open ball in $X$ of center
$x$ and radius $\varepsilon.$ Thus, in particular, we have that
$\operatorname{dist}_{q}(C(X),C(X_{\varepsilon}))\leq\varepsilon$. Let
$n=\#X_{\varepsilon}$. There exists, by Proposition (\ref{easycv}), a Lip-norm
$L_{n}$ on $M_{n}$ such that $\operatorname{dist}_{q}(C(X_{\varepsilon}%
),M_{n})\leq\varepsilon$. By the triangle inequality, we conclude the proof
our proposition. \qed \end{pf}

\bibliographystyle{amsplain}
\bibliography{thesis}

\providecommand{\bysame}{\leavevmode\hbox to3em{\hrulefill}\thinspace}
\providecommand{\MR}{\relax\ifhmode\unskip\space\fi MR }
\providecommand{\MRhref}[2]{%
  \href{http://www.ams.org/mathscinet-getitem?mr=#1}{#2}
}
\providecommand{\href}[2]{#2}
\begin{thebibliography}{10}

\bibitem{BourbakiIntegration}
{N}. {B}ourbaki, \emph{Int\'egration}, {\'E}l\'ements de Math\'ematiques, vol.
  Livre VI, Hermann, Paris, 1963.

\bibitem{Connes89}
A.~Connes, \emph{Compact metric spaces, fredholm modules and hyperfiniteness},
  Ergodic Theory and Dynamical Systems \textbf{9} (1989), no.~2, 207--220.

\bibitem{Connes}
A.~{C}onnes, \emph{Noncommutative geometry}, Academic Press, San Diego, 1994.

\bibitem{Connes97}
A.~Connes, M.~Douglas, and A.~Schwarz, \emph{Noncommutative geometry and matrix
  theory: Compactification on tori}, JHEP \textbf{9802} (1998), hep-th/9711162.

\bibitem{Davidson}
{K}.~{R}. {D}avidson, \emph{C*--algebras by examples}, Fields Institute
  Monographs, American Mathematical Society, 1996.

\bibitem{Dixmier}
J.~{D}ixmier, \emph{Les {C*--}algebres et leur repr\'esentations},
  Gauthier-Villars, 1969, (reprint) Editions Jacques Gabay, 1996.

\bibitem{Dudley}
R.~M. Dudley, \emph{Real analysis and probability}, 2002 ed., Cambridge Studies
  in Advanced Mathematics, vol.~74, Cambridge University Press, 2002.

\bibitem{Glimm62}
{J}. {G}limm, \emph{Families of induced representations}, Pacific J. of Math.
  \textbf{12} (1962), 885--911.

\bibitem{Jaffe}
J.~Glimm and A.~Jaffe, \emph{Quantum physics}, 2nd ed., Springer-Verlag, 1987.

\bibitem{Gromov}
M.~{G}romov, \emph{Metric structures for {R}iemannian and non-{R}iemannian
  spaces}, Progress in Mathematics, Birkh\"auser, 1999.

\bibitem{Hewitt79}
{E}. Hewitt and {K}.~{A}. {R}oss, \emph{Abstract harmonic analysis}, second
  ed., vol.~1, Grundlehren der Mathematischen Wissenschaften, no. 115,
  Springer-Verlag, 1979.

\bibitem{Higson01}
N.~{H}igson, V.~{L}afforgue, and G.~{S}kandalis, \emph{Counter-examples to the
  {B}aum-{C}onnes conjecture},  (2001).

\bibitem{HilSkan87}
M.~{H}ilsum and G.~{S}kandalis, \emph{Morphismes {K}-orient\'es {d'e}spaces de
  feuilles et fonctorialite en th\'eorie de {K}asparov}, Ann. Scient. Ec. Norm.
  Sup \textbf{20} (1987), 325--390.

\bibitem{Kleppner65}
A.~Kleppner, \emph{Multipliers on {A}belian groups}, Mathematishen Annalen
  \textbf{158} (1965), 11--34.

\bibitem{LanRam00}
N.P. {L}andsman and B.~{R}amazan, \emph{Quantization of {P}oisson algebras
  associated to {L}ie algebroids},  (2000), math-ph/0001005.

\bibitem{Latremoliere04}
{F}. {L}atr{\'e}moli{\`e}re, \emph{Finite dimensional approximations of the
  quantum tori for the quantum gromov-hausdorff distance}, Ph.D. thesis,
  University of California at Berkeley, 2004, Advisor: Marc A. Rieffel.

\bibitem{Madore}
J.~{M}adore, \emph{An introduction to noncommutative differential geometry and
  its physical applications}, 2nd ed., London Mathematical Society Lecture
  Notes, vol. 257, Cambridge University Press, 1999.

\bibitem{Zumino98}
B.~Morariu and B.~Zumino, \emph{Super yang-mills on the noncomutative torus},
  Arnowitt Festschrift \textbf{Relativity, Particle Physics, and Cosmology}
  (1998), hep-th/9807198.

\bibitem{PimVoi80a}
{M}. {P}imsner and {D}.~{V}. Voiculescu, \emph{Imbedding the irrational
  rotation algebras into an {AF} algebra}, Journal of Operator Theory
  \textbf{4} (1980), 201--210.

\bibitem{Ram98}
B.~{R}amazan, \emph{Quantification par deformation des vari\'et\'es de
  {L}ie-{P}oisson}, Ph.D. thesis, Universit\'e d'Orleans, France, 1998.

\bibitem{Renault80}
J.~{R}enault, \emph{A groupoid approach to {C*--}algebras}, Lecture Notes in
  Mathematics, vol. 793, Springer-Verlag, 1980.

\bibitem{Renault87}
\bysame, \emph{Repr\'esentation des produits crois\'es d'alg\`ebres de
  groupo\"ides}, J. Operator Theory \textbf{18} (1987), 67--97.

\bibitem{Rieffel89}
M.~A. {R}ieffel, \emph{Continuous fields of {C*--}algebras coming from group
  cocycles and actions}, Math. Ann. \textbf{283} (1989), 631--643.

\bibitem{Rieffel98a}
\bysame, \emph{Metrics on states from actions of compact groups}, Documenta
  Mathematica \textbf{3} (1998), 215--229, math.OA/9807084.

\bibitem{Rieffel99}
\bysame, \emph{Metrics on state spaces}, Documenta Math. \textbf{4} (1999),
  559--600, math.OA/9906151.

\bibitem{Rieffel01}
\bysame, \emph{Matrix algebras converge to the sphere for quantum
  {G}romov--{H}ausdorff distance}, Mem. Amer. Math. Soc. (2001),
  math.OA/0108005.

\bibitem{Rieffel00}
\bysame, \emph{{G}romov-{H}ausdorff distance for quantum metric spaces}, Mem.
  Amer. Math. Soc. \textbf{168} (March 2004), no.~796, math.OA/0011063.

\bibitem{Rudin91}
{W}. {R}udin, \emph{Functional analysis}, 2nd edition ed., McGraw-Hill, 1991.

\bibitem{Szabo01}
R.~Szabo, \emph{Quantum field theory on noncommutative spaces},  (2001),
  hep-th/0109162.

\bibitem{Zeller-Meier68}
{G}. {Z}eller {M}eier, \emph{Produits crois\'es {d'u}ne {C*-}alg\`ebre par un
  groupe {d' A}utomorphismes}, {J}. {M}ath. pures et appl. \textbf{47} (1968),
  no.~2, 101--239.

\end{thebibliography}

\end{document}